\documentclass[10pt,amssymb]{article}
\usepackage{amsfonts,latexsym}
 
\newcommand{\HH}{\mathbb H} 
\newcommand{\R}{\mathbb R} 
\newcommand{\Z}{\mathbb Z} 
\newcommand{\N}{\mathbb N} 
  
\newcommand{\del}{\partial}
\newcommand{\e}{\varepsilon}

\def\ds{\displaystyle}

\newtheorem{theorem}{Theorem} 
\newtheorem{lemma}{Lemma} 
\newtheorem{proposition}{Proposition} 
 
\newtheorem{definition}{Definition} 
\newtheorem{remark}{Remark} 

\begin{document} 
 
\title{Attaching handles to Bryant surfaces}
 
\author{Frank Pacard and Fernando A. A. Pimentel}

\maketitle 

\section{Introduction}

Complete minimal surfaces and constant mean curvature surfaces in Euclidean or hyperbolic space have been the object of much attention in recent years and particular efforts have been devoted in understanding the moduli space of such surfaces. More precisely, the surfaces we are interested in are complete, noncompact surfaces which may or may not be embedded but enjoy the somehow weaker, and slightly different, property to be Alexandrov embedded (since we do not need this notion later on in the paper, we simply refer to  \cite{Cos-Ros} or \cite{Gro-Kus-Sul} for a precise definition of Alexandrov embeddedness). 

\medskip

In the case of minimal surfaces of the Euclidean 3-space ${\R}^3$, a fairly good picture of the space of complete Alexandrov embedded minimal surfaces with finite total curvature is now available. We will denote by $g$ the genus of such a surface and by $n$ its number of ends, which are all known to be of catenoidal type \cite{Cos-Ros}. To begin with,  we can use the Weierstrass representation theory to obtain the existence of these surfaces for any number of ends $n \geq 2$, when the genus is $0$ \cite{Jor-Mee}, \cite{Kat-Ume-Yam}, and also to obtain some examples for higher genus \cite{Ber-Ros}. Let us also mention the recent result of S.D. Yang which provides a wealth of new examples be settling a connected sum theory for these surfaces. This construction allows one to add as many ends as one wants to a given minimal surface and also to connect sum any two of such surfaces. However, and to our knowledge, it is not known whether such minimal surfaces always exist for any number of  ends $n \geq 2$ and any genus $g\geq 1$, though these surfaces are conjectured to exist.  C.P. Cosin and A. Ros have completed the classification of these minimal surfaces when all the ends have their axis in the same plane and when the genus is $0$. Finally, let us mention that J. Perez and A. Ros have determined the structure of the moduli space of these surfaces \cite{Per-Ros}.  

\medskip

In the case of  complete Alexandrov embedded constant mean curvature surfaces in the Euclidean 3-space ${\R}^3$,  beside the well known one parameter family of surfaces obtained by Delaunay \cite{Del} in the middle of the nineteenth century,  examples have been obtained by N. Kapouleas \cite{Kap}, K. Grosse-Brauckmann \cite{KGB} and R. Mazzeo and F. Pacard \cite{Maz-Pac}. Let us also mention the result of R. Mazzeo, F. Pacard, D. Pollack and J. Ratzkin \cite{Maz-Pac-Pol-Rat}  which parallel, for constant mean curvature surfaces, the result of S.D. Yang for minimal surfaces.  As above this result  allows one to add as many ends as one wants to a given constant mean curvature surface and also to connect sum two of these surfaces. The moduli space theory was developed by R. Kusner, R. Mazzeo and D. Pollack \cite{Kus-Maz-Pol} and the classification of complete constant mean curvature surfaces, of genus $0$, with $3$ ends was recently completed by K. Grosse-Brauckman, R. Kusner and J. Sullivan \cite{Gro-Kus-Sul}.

\medskip

Finally, we turn to the main subject of the present paper : complete embedded constant mean curvature-$1$ surfaces in the hyperbolic $3$-space ${\HH}^3$.  A representation formula for constant mean curvature-$1$ surfaces ({\small CMC}-1 for short) in hyperbolic $3$-space has been discovered by R.L. Bryant \cite{Bry}. This shows that, to some extent, {\small CMC}-$1$ surfaces in ${\HH}^3$ behave like minimal surfaces in ${\R}^3$.  Beside the well known horospheres, which have one end, there exits a one parameter family of {\small CMC}-$1$ surfaces of revolution with $2$ ends known as "catenoid cousins" \cite{Bry}.  Using the Weierstrass representation of {\small CMC}-$1$ surfaces, obtained by R.L. Bryant, a number of geometrically interesting surfaces (eventually immersed) have been found. We refer to the recent paper \cite{Ros-Ume-Yam-2} for further informations about these examples. 

\medskip

In the present paper, we would like to give an answer to some natural questions about complete embedded {\small CMC}-$1$ surfaces in hyperbolic space. Before we do so let us introduce the following~:
\begin{definition}
A constant mean curvature-$1$ surface in ${\mathbb H}^3$ is said to have regular ends if its hyperbolic Gauss map extends through the punctures.
\end{definition}
The precise definition of the hyperbolic Gauss map of a surface in hyperbolic space is given in \cite{Bry}. All the surfaces we will be considering in this paper do have regular ends since their ends will always be modeled after the end of a catenoid cousin. 
\begin{definition}
Given $g \geq 0$ and $n \geq 1$, the unmarked moduli space of {\small CMC}-$1$ surfaces ${\cal M}_{g,n}^u$ in ${\HH}^3$ is defined to be the moduli space of complete, embedded, constant mean curvature-$1$ surfaces of genus $g$  in ${\mathbb H}^3$, which are connected and have $n$ regular ends.
\end{definition}

\medskip

When $g=0$, ${\cal M}^u_{0,n}$ is fairly well understood. For example, when $n=1$, it consists of only one point : the horosphere, when $n=2$, it consists in the one parameter family of catenoid cousins  which have already been mentioned (see \cite{Lev-Ros} and \cite{Ume-Yam}). Finally, when $n=3$, the recent work of P. Collin and H. Rosenberg shows that ${\cal M}^u_{0,3}$  is homeomorphic to the moduli space of geodesic triangles in ${\R}P^3$ \cite{Col-Ros}, paralleling what is already known for minimal and constant mean curvature surfaces in ${\R}^3$. Beside this result, let us mention the result of  W. Rossman, M. Umehada and K. Yamada \cite{Ros-Ume-Yam}, showing for example that ${\cal M}^u_{1, 3}$ is not empty. Their proof relies on the existence of a genus $1$ minimal surface with $3$ ends in ${\R}^3$ which is non degenerate, and  it does not extend easily to arbitrary genus. Granted the above results, the first question we would like to address is~:

\medskip

For which values of $g$ and $n$, is the moduli space ${\cal M}^u_{g,n}$ not empty ? 

\medskip

\noindent
Regarding this problem, we have obtained the
\begin{theorem}
The moduli space ${\cal M}_{g , n}^u$ is not empty in the following cases~: 
\begin{enumerate}
\item[(i)]  When $g=0$ and  $n \geq 1$.
\item[(ii)] When $g \geq 1$ and  $2 \, n \geq   g+5$.
\end{enumerate}
\label{th:1}
\end{theorem}
As will become clear in the remaining of the paper, this result, as well as the next one will be obtained by desingularizing a finite number of horospheres. We will say that the horospheres ${\cal H}^1, \ldots, , {\cal H}^n$ are "either disjoint or tangent" if any two horospheres ${\cal H}^i$ and ${\cal H}^j$ are either disjoint or their intersection reduces to the point where they are tangent. 

\medskip

In order to explain how the constraint (ii) is obtained, we define $n_g$ to be the minimal number of either disjoint or tangent horospheres which are needed for their connected sum to be an embedded surface of genus $g$. It is easy to see that, in order to obtain a genus $1$ surface, at least three horospheres ${\cal H}^1, {\cal H}^2 , {\cal H}^3$, each of which is tangent to the other two, are necessary. Now, one can slightly reduce the radii of the three horospheres ${\cal H}^1, {\cal H}^2 , {\cal H}^3$ and consider their connected sum at the three points where they were tangent to produce an embedded genus $1$ surface. Hence, we already have $n_1 \geq 3$. Next, we assume that we have already found a configuration of $n$ either disjoint or tangent horospheres ${\cal H}^1, \ldots, {\cal H}^n$ whose connected sum produces an embedded genus $g$ surface (observe that it is not always necessary to connect the horospheres at every point where they are tangent to obtain a surface of genus $g$, and, since we slightly reduce the radii of the horospheres, the connected sum surface will not be singular regardless of the fact that we have connected the surfaces at every point or not). Now, one can find an horosphere ${\cal H}^{n+1}$ which is tangent to (at least) three of the initial configuration of horospheres ${\cal H}^1, \ldots, {\cal H}^n$. The connected sum of ${\cal H}^1, \ldots, {\cal H}^{n+1}$ obtained by slightly reducing the radii, connecting ${\cal H}^1, \ldots, {\cal H}^n$ as before and then connecting the additional horosphere ${\cal H}^{n+1}$, whose radius is slightly reduced too, at two points (resp. three points) produces an embedded genus $g+1$ (resp. $g+2$) surface.  Hence we obtain the inequalities $n_{g+1}\geq n_g+1$ and $n_{g+2}\geq n_g +1$. The inequality in (ii) follows at once by induction.

\medskip

Paralleling what has been done for minimal or constant mean curvature of ${\R}^3$, we can prove the~:
\begin{theorem}
For any $g\geq 0$ and any $n\geq 1$, the moduli space ${\cal M}^u_{g,n}$ is a real analytic variety of formal dimension $3\, n$ (before moding out by the action of isometries of ${\HH}^3$).
\label{th:2}
\end{theorem}
This formal dimension is achieved at any unmarked nondegenerate point of the moduli space. In addition to the above result, we prove that all the surfaces constructed in this paper are marked nondegenerate provided $g=0$.

\medskip

Before stating our next result, we give the~:
\begin{definition}
For any $g \geq 0$ and any $(p_i)_{1\leq i \leq n}$, finite set of distinct points of $\del {\mathbb H}^3$, the marked moduli space of Bryant surfaces of genus $g$  is denoted by ${\cal M}_{g, (p_i)_i}^m$ and defined to be the moduli space of embedded {\small CMC}-$1$ surfaces of genus $g$ in ${\mathbb H}^3$ which are complete, connected and have ends at the points $p_i$.
\end{definition}
Granted the above definition, a natural question is~: 

\medskip

For which sets of points $(p_i)_{1\leq i\leq n} \in \del {\mathbb H}^3$ and for which genus $g$ is the moduli space ${\cal M}_{g , (p_i)_i}^m$ not empty ?

\medskip

This question can also be understood as a Plateau problem at infinity for {\small CMC}-$1$ surfaces in ${\R}^3$ since it amounts to look for a surface with prescribed asymptotic behavior at infinity. As a byproduct of our construction, we obtain the~: 
\begin{theorem}
Given any  $(p_i)_{1\leq i\leq n}$ finite set of points of $\del {\mathbb H}^3$,  the moduli spaces ${\cal M}^m_{0, (p_i)_i}$ and  ${\cal M}^m_{1, (p_i)_i}$ are not empty if $n \geq 3$.
\label{th:3}
\end{theorem}
Again the result will be obtained by desingularizing a finite number of horospheres which are either tangent or disjoint. Given the points $p_i$, it is easy to check that we can find horospheres with ends at the points $p_i$, in such a way that their connected sum is an embedded surface of genus greater than or equal to $1$. In fact, one can find $3$ horospheres ${\cal H}^1, \ldots, {\cal H}^3$ with ends at $p_1, \ldots, p_3$, each of which is tangent to the other two. Now, one can slightly reduce the radii of the three horospheres ${\cal H}^1, \ldots, {\cal H}^3$ and consider their connected sum at two or three of the three points where ${\cal H}^1, \ldots, {\cal H}^3$ are tangent to produce an embedded genus $0$, or genus $1$, surface. Next, we assume that we have already found a configuration of $n$ either disjoint or tangent horospheres ${\cal H}^1, \ldots, {\cal H}^n$, with ends at $p_1, \ldots, p_n$ respectively,  whose connected sum produces an embedded genus $0$, or of genus $1$, surface (again, observe that to produce an  embedded genus $0$, or a genus $1$, surface it is usually not necessary to connect the horospheres at every point where they are tangent). One can find an horosphere ${\cal H}^{n+1}$, with end at $p_{n+1}$ which is tangent to (at least) one of the initial configuration of horospheres ${\cal H}^1, \ldots, {\cal H}^n$. The connected sum of ${\cal H}^1, \ldots, {\cal H}^{n+1}$ obtained by connecting ${\cal H}^1, \ldots, {\cal H}^n$ as before and connecting the additional horosphere ${\cal H}^{n+1}$ at one point produces an embedded genus $0$, or genus $1$, surface. 

\medskip

This result extends to any countable set of points but we shall not insist on this extension.  We have~: 
\begin{theorem}
Given $g \geq 0$ and given any finite set of points $(p_i)_{1 \leq i \leq n}$ in $\del {\mathbb H}^3$, the moduli space ${\cal M}^m_{g, (p_i)_i}$ is a real analytic variety of formal dimension $n$ (before moding out by the action of isometries of ${\HH}^3$).
\label{th:4}
\end{theorem}
Again, this formal dimension is achieved at any marked nondegenerate point of the moduli space. We prove that all the surfaces constructed in this paper are marked nondegenerate provided $g=0$ and for a generic choice of the points $(p_i)_{1\leq i \leq n}$.

\medskip

\noindent
{\bf Acknowledgment~:} Part of this paper was written while first author was visiting the Mathematical Sciences Research Institute in Berkeley. He would like to take this opportunity to thank the MSRI for their support and hospitality. This paper was written while the second author held a post-doc position in the Universit\' e Paris
 XII. He would like to thank the mathematical department of Paris XII for his hospitality. The authors would also like to thank  H. Rosenberg for bringing this problem to their attention.

\section{Preliminaries}

\subsection{Usual models for ${\HH}^3$}

We recall the $3$ usual models which are used to describe ${\HH}^3$. More information can be found in \cite{Spi} or in $\cite{Bry}$.

\medskip

\noindent
{\bf The ball model}. In this model the hyperbolic space is  described by the unit open ball of ${\R}^3$. So
\[
{\HH}^3 :=\{ x\in {\R}^3 \quad : \quad |x| <1\},
\]
endowed with the metric $g_{hyp} : = 4 \, (1-|x|^2)^{-2} \, dx^2$. While this model is probably the best one for visualizing pictures of {\small CMC}-$1$ surfaces in ${\HH}^3$ (see for example \cite{Ros-Sat}), computations are not so easy to perform in it.

\medskip

\noindent
{\bf The Minkowski model}. In this model the hyperbolic space is described by one sheet of a two sheeted hyperbola in the Lorentz space ${\mathbb L}^4$. More precisely, we have 
\[
{\HH}^3 :=  \{ (x_0, x) \in {\R}^{1+3}\quad : \quad x_0^2 = 1+|x|^2, \quad x_0>0 \},
\]
endowed with the induced metric of ${\mathbb L}^4$ given by $g_{lor}:= -dx_0^2+dx^2$. This model is particularly useful when working with the Weierstrass representation \cite{Bry}. For example, in this model the one parameter family of catenoid cousins can be parameterized by
\begin{equation}
\begin{array}{rllll}
x_0 (s, \theta)  & = & \ds \frac{1}{4 (t-1)} \, \left( t^2 \, \cosh ((t-2) s ) + (t-2)^2 \, \cosh (t s) \right)\\[3mm]
x_1 (s, \theta)  & = & \ds \frac{1}{2 (t-1)} \, t \, (t-2) \,  \cosh ((t-1)s)  \, \cos \theta\\[3mm]
x_2 (s, \theta)  & = & \ds \frac{1 }{2(t-1)} \,  t \, (t-2) \,  \cosh ((t-1)s)  \, \sin \theta\\[3mm]
x_3 (s, \theta)  & = & \ds \frac{1}{4 (t-1)} \, \left( t^2 \, \sinh ((t-2)s) - (t-2)^2 \, \sinh (t s) \right) 
\end{array}
\label{eq:2-0}
\end{equation}
for any $t >1$, $t \neq 2$.  It turns out that, when $t \in (1,2)$ the catenoid cousins are embedded and, when $t$ tends to $2$, this one parameter family converges to the union of two horospheres, while, when $t$ tends to $1$, the surfaces tend to $\del {\HH}^3$ uniformly. This family continues for $t >2$ but the surfaces are not embedded anymore. This phenomena is very reminiscent to what happens in Euclidean space for {\small CMC} surfaces. For example, Delaunay surfaces come in a one parameter family for some parameter $\tau \leq 1$. The embedded members of this family are called unduloids and interpolate between a cylinder, when $\tau =1$, to a bead of spheres arranged along an axis, when $\tau=0$. Beyond this limit, when $\tau <0$, this family continues but the surfaces, which are called nodoids, are not embedded anymore \cite{Maz-Pac}. 

\medskip

\noindent
{\bf The Poincar{\'e} or upper half space model}. In this last model, the hyperbolic space is identified with ${\mathbb H}^3 := \{(x,z) \in {\R}^2\times {\R} \quad : \quad z >0\}$, endowed with the metric
\[
g_{hyp} := \frac{1}{z^2}\, (dx^2 + dz^2).
\]
We can also identify $\del {\HH}^3$ with ${\R}^2$. This is this last model which will be most useful for us. Hence, in all the remaining of the paper, we will work with it. 

\subsection{Isometries of ${\mathbb H}^3$}

We recall the theorem of Liouville concerning the classification of isometries of ${\mathbb H}^3 := \{(x,z) \in {\R}^2\times {\R} \quad : \quad z >0\}$, endowed with the metric
\[
g_{hyp} := \frac{1}{z^2}\, (dx^2 + dz^2).
\]

\begin{theorem}
The isometries of $({\mathbb H}^3, g_{hyp})$ are the restrictions to ${\mathbb H}^3$ of conformal transformations of ${\R}^3$ that take ${\mathbb H}^3$ onto itself.
\end{theorem}
Hence, the isometries of ${\mathbb H}^3$ are compositions of horizontal translations
\[
{\cal T}_a (x,z) :=  (x+a, z),
\]
for any $a \in {\R}^2$, dilations centered at the origin
\[
{\cal D}_\lambda (x,z) : =  \lambda \, (x,z),
\]
for any $\lambda >0$, elements of the orthogonal group
\[
{\cal R}_{R} (x,z) := (R \, x, z),
\]
for any $R \in O(2)$ and the inversion centered at the origin 
\[
{\cal I}_0 (x,z) : = \frac{(x,z)}{|x|^2+z^2}.
\] 

\subsection{Mean curvature in hyperbolic space}

We still consider the upper half space model ${\mathbb H}^3 = \{ (x,z)  \in {\R}^2 \times {\R} \, : \, z >0\}$, endowed with the metric 
\[
g_{hyp} : = \frac{1}{z^2}\, ( dx^2 +  dz^2).
\]
In this model, the mean curvature $H_{hyp}$ of a surface $\Sigma$, endowed with the metric induced by $g_{hyp}$, can be compared to the mean curvature $H_{eucl}$ of the same surface which this time is considered to be in ${\R}^3$, and hence is endowed with the metric induced by $g_{eucl}$. This is the content of the following classical result whose proof can be found in \cite{Tou}~:
\begin{proposition}
Let $\Sigma$ be a surface contained in the upper half space. We denote by $z$ the height function. If $N_{eucl} $denotes the normal to the surface in $({\R}^3, g_{eucl})$, we denote by $N^z_{eucl}$ the coordinate of $N_{eucl}$ along the $z$ axis. Then the mean curvatures of $\Sigma$ in $({\mathbb H}^3, g_{hyp})$ and in $({\R}^3, g_{eucl})$ are related by
\begin{equation}
H_{hyp} = z \, H_{eucl} + N^z_{eucl} .
\label{eq:2-1}
\end{equation}
\end{proposition}
{\bf Proof :} If $x_1, x_2$ and $z$ are the standard coordinates in ${\mathbb H}^3$, the tangent space is spanned by the vectors
\[
\del_1 := \del_{x_1}, \qquad \del_2 : = \del_{x_2} \qquad \mbox{and}\qquad \del_3 : = \del_z.
\]
The Christoffel symbols associated to $g_{hyp}$ all vanish except
\[
\Gamma_{11}^3 = \Gamma_{22}^3 = \frac{1}{z},
\] 
and 
\[
\Gamma_{13}^1 = \Gamma_{31}^1 = \Gamma_{23}^2= \Gamma^2_{32} = \Gamma_{33}^3 = - \frac{1}{z}.
\] 
Now, if $X  : = \sum_i X^i \, \del_i$ and if $Y : = \sum_j Y^j \, \del_j$ are two tangent vector fields, we have the expression of the covariant derivative in $({\mathbb H}^3, g_{hyp})$ which is given by
\[
\overline\nabla_X Y =  \sum_k \, \left( \sum_{i,j} \, X^i \, Y^j \, \Gamma^k_{ij} + X (Y^k) \right)\, \del_k .
\]
We find explicitely
\[
\overline \nabla_X Y  = \sum_k X (Y^k) \, \del_k - \frac{Y^3}{z} \, X - \frac{X^3}{z} \, Y + \frac{1}{z}\, \left( \sum_i X^i \, Y^i \right)  \, \del_3 .
\]
Hence, if $Y$ is a vector field normal to $\Sigma$ and if $X$ is a vector field tangent to $\Sigma$, we find that 
\[
[ \overline \nabla_X \, Y]^{T} =  \left[ \sum_k X (Y^k) \, \del_k \right]^T  - \frac{Y^3}{z} \, X,
\]
where $[\, \cdot \, ]^T$ denotes the orthogonal projection onto $T \, \Sigma$.

\medskip

In the particular case where $Y$ is the normal huperbolic vector field $N_{hyp}$, we obtain
\[
[ \overline \nabla_X \, N_{hyp} ]^{T} =  [ \nabla_X  \, N_{hyp}]^T  - \frac{N^z_{hyp}}{z}\, X,
\]
where $\nabla_X Y$ is the covariant derivative of $Y$ along $X$ in $({\mathbb R}^3, g_{eucl})$ and where $N^z_{hyp}$ is the coordinate along the $z$ axis of the normal hyperbolic vector $N_{hyp}$. Finally, we use the fact that the normal vectors $N_{eucl}$ and $N_{hyp}$ are related by the identity
\[
N_{hyp} = z \, N_{eucl},
\]
to conclude that
\[
\begin{array}{rllll}
[\overline \nabla_X N_{hyp}] & = & [\nabla_X (z \, N_{eucl})]^T - N_{eucl}^z \, X \\[3mm]
			         & = & [ z \, \nabla_X N_{eucl} + X(z) \, N_{eucl}]^T - N_{eucl}^z \, X \\[3mm]
			         & = & z \, [\nabla_X N_{eucl}]^T -N^z_{eucl} \, X.
\end{array}
\]
The mean curvature in $({\mathbb H}^3, g_{hyp})$ is defined to be half the trace of the mapping
\[
X \longrightarrow - \, [\overline \nabla_X N_{hyp}]^T,
\]
while the mean curvature in $({\mathbb R}^3, g_{eucl})$ is defined to be half the trace of the mapping
\[
X \longrightarrow - \, [\nabla_X N_{eucl}]^T. 
\]
Ganted these defintions, the relevant formula follows at once. \hfill $\Box$

\medskip

In particular, the mean curvature of the graph of a positive function $u$ in $({\mathbb H}^3, g_{hyp})$ is given by
\begin{equation}
H (u) :=  \frac{1}{\sqrt{1+ |\nabla u|^2}}+ \frac{u}{2}\, \mbox{div}\left( \frac{\nabla u}{\sqrt{1+ |\nabla u|^2}}\right)
\label{eq:2-2}
\end{equation}
Observe, when the normal vector is assumed to point upward,  that the plane $z=z_0$  has constant mean curvature equal to $1$  and so has the lower half of the horosphere
\[
z = R - \sqrt{R^2 -x^2-y^2}.
\]
While, because of the chosen orientation, the upper half of the horosphere 
\[
z = R + \sqrt{R^2 -x^2-y^2},
\]
has mean curvature $-1$.

\subsection{The ends of a catenoid cousin}

Interesting for us will be a classification of the behavior of a catenoid cousin near one of its ends, since it is known, from the work of P. Collin, L. Hauswirth and H. Rosenberg \cite{Col-Hau-Ros}, that they model the ends of the surfaces we are interested in.  Obviously, an explicite description of these ends can be obtained by considering the expression of the catenoid cousins which is given in (\ref{eq:2-0}) for the Minkowski model and pulling this back to the upper half space model. Nevertheless, we prefer to provide a shorter  indirect proof which amounts to look for rotationally symmetric graphs in the upper half space model for ${\HH}^3$.  Hence, we want to look for radially symmetric solutions of the equation  $H(u)=1$, where $H(u)$ is given in (\ref{eq:2-2}). The idea is to fix $t >1$ and to look for a solution $u$ of the form
\[
u(x) : = r^{2-t} \, (1+v (r)),
\]
where the function $v$ is small (compared to $1$) and only depends on $r$.  Since the functions we are looking for are radial, the nonlinear partial differential equation $H(u)=1$ reduces to an ordinary differential equation which we can solve by using a standard fixed point procedure. As a result, we obtain that, given $1 <t_0 < t^0 <\infty$, there exists  $r_0 >0$  and, for all $t \in [t_0, t^0]$, there exists a solution of (\ref{eq:2-2}) of the form
\[
u_t : = r^{2-t} \, (1+ v_t)
\]
which is defined over ${\R}^2 - B_{r_0}$, where the function $v_t$ satisfies for all $k\geq 0$ 
\[
 |\nabla^k v_t |\leq  c_k \, r^{2-2t-k},
\]
for some constant $c_k>0$ independent of $t\in [t_0, t^0]$. This way we obtain, up to a dilation factor,  a local parameterization of the end of a (rotationally symmetric) catenoid cousin $C_t$ in the upper half space model. In particular, when $t=2$, the function $u_2 \equiv 1$ parameterizes the end of a horosphere.

\medskip

Let us denote by $L_t$ the Jacobi operator, that is the linearized mean curvature operator about $C_t$ with respect to the hyperbolic normal vector field. It is well known that
\[
L_t : = \Delta_{C_t} + |A_{C_t}|^2 + \mbox{Ric}_{{\HH}^3} (N_{hyp} , N_{hyp}),
\]
where $|A_{C_t}|^2$ denotes the norm of the second fundamental form, $\mbox{Ric}_{{\HH}^3}$ is the Ricci tensor of $({\HH}^3, g_{hyp})$ and $N_{hyp}$ is a unit hyperbolic normal  vector to $C_t$. It will be important, in the study of the moduli spaces, to understand the asymptotic behavior of some Jacobi fields, i.e. solutions of the homogeneous problem $L_t w = 0$. These Jacobi fields correspond to explicit one-parameter geometric transformations of the catenoid cousin $C_t$, say $\zeta \rightarrow C_t(\zeta)$ with $C_t(0) = C_t$. For all $\zeta$ small enough, $C_t(\zeta)$ can be written (at least locally) as a normal graph over $C_t$ and differentiation with respect to $\zeta$ gives rise to one Jacobi field which is globally defined on $C_t$, however, we will only be interested by its behavior near one end of $C_t$, for example, near $\infty$.

\medskip

Using the above procedure, if one considers $t' \rightarrow C_{t+t'}$ as a one parameter family of transformation of a fixed catenoid cousin $C_t$, one finds the Jacobi field $h^{0,+}$, whose behavior near $\infty$ is given by
\begin{equation}
h^{0,+} \sim \log r ,
\label{eq:2-5}
\end{equation}
while  dilations induce the one parameter family $\lambda \longrightarrow {\cal D}_\lambda \, C_t$ and  give rise to the Jacobi field $h^{0,-}$, whose behavior near $\infty$ is given by 
\begin{equation}
h^{0,-} \sim 1 .
\label{eq:2-6}
\end{equation}
Next, translations induce the family  $a \longrightarrow {\cal T}_a \, C_t$  and give rise to the linearly independent Jacobi fields $h^{-1,\pm}$, whose behavior near $\infty$ is given by 
\begin{equation}
h^{-1, \pm} \sim r^{-1}\, e^{\pm i \theta} .
\label{eq:2-7}
\end{equation}
Finally using a composition of translations and an inversion, we find the family $a \longrightarrow {\cal T}_{\tilde a} \circ {\cal I}_0 \circ {\cal T}_a \,  C_t$, where $\tilde a := -\frac{a}{|a|^2}$  and give rise to the two other independent Jacobi fields $h^{+1,\pm}$, whose behavior near $\infty$ is given by 
\begin{equation}
h^{+1, \pm} \sim r \, e^{\pm i \theta} .
\label{eq:2-8}
\end{equation}
The key observation is that, up to higher order terms, all these Jacobi fields behave like harmonic functions. Though this is not explicite in the notation, all these Jacobi fields do depend on $t$.

\section{Basic preparation for the construction}

We consider a family of horospheres $({\cal H}^i)_{i\in I}$ which are
either tangent or disjoint. Recall that either ${\cal H}^i$ and ${\cal H}^j$ are disjoint or their intersection reduces to the point where they are tangent.  Assume that $\cup_{i\in I} \, {\cal H}^i$ is a connected subset of ${\mathbb H}^3$. We will denote by $r_i$ the radius of the horosphere ${\cal H}^i$ and by $a^i_\infty \in \del {\mathbb H}^3$, the end of the horosphere. Possibly ${\cal H}^i$ is the horosphere given by the equation $z=\mbox{cte}$, in which case we will have $a_\infty^i =\infty$.

\medskip

We will denote by $J =\{(i,j) \, : \, i \neq j , \quad  {\cal H}^i \cap {\cal H}^j \neq \emptyset\}$ the set of couples $(i,j)\in I\times I$ such that the horospheres ${\cal H}^i$ and ${\cal H}^j$ are tangent, in which case, the point of tangency between ${\cal H}^i$ and ${\cal H}^j$ will be denoted $p_{ij}$. In particular, $p_{ij} = p_{ji}$. Now, in the forthcoming gluing procedure, we may wish to connect the two horospheres ${\cal H}^i$ and ${\cal H}^j$ at the point $p_{ij}$ or we may not. The point being that we have the freedom to choose a subset of indices $\hat J \subset J$ which correspond to the points $p_{ij}$ where we will actually connect ${\cal H}^i$ and ${\cal H}^j$. The subset $\hat J$ can be arbitrary, the only restriction being that we want the collection of horospheres $({\cal H}^i)_{i\in I}$ connected at the points $p_{ij}$, for all $(i,j) \in \hat J$ to be a connected set. Finally, given a horosphere ${\cal H}^i$, we will denote by $I_i$ the set of indices $j \in I$ such that  the horosphere ${\cal H}^j$ has to be connected to  ${\cal H}^i$, in other words for which $(i,j) \in \hat J$.
 
\medskip

{\bf Definition of the parameters $a^i_j$ and $d^i_j$~:} We can now introduce the parameters we will need. For any $\e\in (0,1)$, we define $\eta >0$ by the formula
\begin{equation}
\eta  := - \e \, \log \e .
\label{eq:3-1}
\end{equation}
Now, for each horosphere ${\cal H}^i$, one can perform an isometry ${\cal J}_i$ so that the image of ${\cal H}^i$ by ${\cal J}_i$ becomes the horosphere of radius $\frac{1}{2}$, with a end at the origin. Granted this, we define
\[
{\cal H}^i (\e) : = {\cal J}_i^{-1}\circ {\cal D}_{\frac{1}{1+\eta}}\circ {\cal J}_i \, ({\cal H}^i ).
\]
What we have actually done, is a reduction of the radii of all the horospheres ${\cal H}^i$ in a consistent way.

\medskip

Since the analysis will be the same for each horosphere,  we will focus our attention on the gluing of the horospheres which have to be connected to ${\cal H}^i$. Performing an isometry, if this is necessary, we can assume that ${\cal H}^i$ is the horosphere given by $z=1$ with end at $a^i_\infty = \infty$ and hence all horospheres ${\cal H}^j$ which have to be connected to ${\cal H}^i$ are horospheres of radius $\frac{1}{2}$ with end at $a^j_\infty \in \del {\HH}^3$. Observe that, in this case, ${\cal H}^i(\e)$ is the plane $z = 1+\eta$ and, for all $j \in I_i$,  ${\cal H}^j(\e)$ is the horosphere or radius $\frac{1}{2(1+\eta)}$ with end at $a^j_\infty$.

\medskip

Given any point $a_{j}^i \in  \del {\HH}^3$ close to $a^i_\infty= \infty$ and any point $a_{i}^j \in \del {\HH}^3$ close to $a^j_\infty$ (we assume that $a_j^i\neq a_i^j$), we will denote by $\Gamma_{ij}$ the geodesic in ${\HH}^3$ with end points $a_j^i$ and $a_i^j$. If $a_j^i$ and $a_i^j$ are respectively close enough to $a^i_\infty$ and $a^j_\infty$, say $|a^i_j -a^j_\infty|>1$ and $|a^j_i-a^j_\infty| <1$, this  geodesic is transverse to the half sphere of radius $1$ centered at the origin,$S_0 := \{(x,z) \in {\HH}^3\, : \, |x|^2+ z^2 =1\}$. Hence we can define the point  $q_{ij}$ to be the intersection between the geodesic $\Gamma_{ij}$ and the half sphere $S_0$.
\[
\{ q_{ij} \} : = S_0 \cap \Gamma_{ij}.
\]
Furthermore, if $\e$ is small enough, we can define $q^i_j (\e)$ (resp. $q_i^j(\e)$) to be the point close to $q_{ij}$ where the geodesic $\Gamma_{ij}$ intersects ${\cal H}^i (\e)$ (resp.  ${\cal H}^j(\e)$). The projection of $q^i_j(\e)$ over $\del {\HH}^3$ will be denoted $x^i_j(\e) \in \del {\HH}^3$.

\medskip

Now, for any point $p \in \Gamma_{ij}$, we will denote by $d_{j}^i (p)$ the signed geodesic distance between $p$ and $q_{ij}$
\[
d^i_{j} (p) :=  \mbox{dist}_{{\HH}^3} (p, q_{ij}).
\]
That is we consider that the geodesic $\Gamma_{ij}$ is oriented and, for example, we assume that $d_j^i (p) \geq 0$  when $p$ belongs to the geodesic arc joining $a^i_j$ to $q_{ij}$ and  $d_j^i  (p) \leq 0$ when $p$ belongs to the geodesic arc joining $q_{ij}$ to $a^j_i$. Observe that, once an orientation of the geodesic $\Gamma_{ij}$ is chosen, we have $d^i_j (p)= -d^j_i(p)$.

\medskip

Hence, for each couple $(i, j) \in \hat {\cal J}$, we have defined so far the parameters $a^i_j, a^j_i \in \del {\HH}^3$ and $d^i_j = -d^j_i\in {\R}$.

\medskip

{\bf Analysis of an isometry ${\cal J}^i_j$~:} We define some isometry which depends on all the parameters we have just defined and analyze precisely its effect in some fixed neighborhood of $q_{ij}$.  To do so, we restrict the parameters so that they satisfy
\begin{equation}
|a^i_j| \geq \frac{1}{\kappa \, \e}, \qquad  |a^j_i-a^j_\infty| \leq \kappa \, \e, \qquad \mbox{and}\qquad  |d^i_j| \leq \kappa \, \e ,
\label{eq:*}
\end{equation}
where $\e \in (0,1)$ and where the constant $\kappa$ will be fixed, independently of $\e$, at the end of \S 6.  The first inequality in (\ref{eq:*}) may seem rather surprising, however, observe that this inequality just ensures that $a^i_j$ is close to $a_\infty^i = \infty$. To convince oneself, one can perform an inversion and obtain an inequality similar to the second inequality in (\ref{eq:*}). If $(x_{ij}, z_{ij})$ are the coordinates of the point $\tilde p_{ij} \in \Gamma_{ij}$, such that $d^i_j ( \tilde p_{ij})=d^i_j$, we define $\lambda^i_j >0$ by the identity
\[
\lambda^i_j \, z_{ij}  = |x_{ij}-a_j^i|^2+ z_{ij}^2 .
\]
For convenience, we set  $\xi^i_j : = a^i_j- a^j_i$ and consider the planar symmetry 
\[
{\cal S}^i_j (x,z) : =  \left( x - \frac{2}{|\xi^i_j|^2} \, (x-a^j_i)\cdot \xi^i_j \, \xi^i_j , z \right) .
\]
Granted these notations,  we now consider the isometry 
\begin{equation}
{\cal J}^i_j  (x,z) : = {\cal S}^i_j \left(  \lambda_j^i \, \left( \frac{(x- a^i_j, z)}{|x-a^i_j|^2 + z^2} - \frac{(a^j_i - a^i_j, 0)}{|a^j_i - a^i_j|^2}\right) + (a^j_i, 0) \right).
\label{eq:ll}
\end{equation}
The mapping ${\cal J}^i_j$ is the unique isometry of ${\HH}^3$ which sends the point $a^i_j$ to $\infty$, keeps the point $a^j_i$ fixed and for which, by construction, the image of the point $p_{ij} = (x_{ij}, z_{ij})$ by ${\cal J}^i_j$ lies in the plane $z=1$. We would like to obtain an expansion of ${\cal J}^i_j$. To this aim, we first observe that, thanks to (\ref{eq:*}), the following expansions hold
\[
z_{ij}  =  1 +  d^i_j + {\cal O}_\kappa (\e^2).
\]
Moreover, we have 
\[
x_{ij}  =  a^j_i + {\cal O}_\kappa  (\e) ,
\]
and hence we obtain the expansion for $\lambda^i_j$
\[
\lambda^i_j =  \frac{|a^i_j -a^j_i|^2}{1+d_j^i } (1 + {\cal O}_\kappa (\e^2)) .
\]
Finally, with little work one finds the asymptotic expansion of ${\cal J}^i_j $
\[
{\cal J}^i_j  (x,z)   = \ds \left( x  + {\cal O}_\kappa (\e ),   z - d^i_j + \frac{2}{|\xi^i_j|^2} \, (x- a^j_i)\cdot  \xi^i_j+ {\cal O}_\kappa ( \e^2 \, \log \e )\right) ,
\]
which is valid for all $(x,z) \in {\HH}^3$ satisfying $|x -a^i_j|\leq 1$ and $|z-1|\leq 2 \, \eta = - 2 \, \e \, \log \e$. Here the subscript $\kappa$ in ${\cal O}_\kappa(\cdot)$ is intended to point out that this quantity may depend on $\kappa$.

\medskip

{\bf Strategy~:} Roughly speaking our strategy can now be described as follows : For each $(i,j) \in \hat J$, we will choose $a^i_j \in \del {\HH}^3$ close to $\infty$, $a^j_i \in \del {\HH}^3$ close to $a^j_\infty$ and $d^i_j \in {\R}$. Then, for each $(i,j) \in \hat J$, we define ${\cal J}_j^i$ as in (\ref{eq:ll}) and we "insert" between the images of the two horospheres ${\cal H}^i(\e)$ and ${\cal H}^j(\e)$ by this isometry, a properly rescaled and "symmetrized"  vertical catenoid which is centered at the point $(a^i_j, 1)\in {\HH}^3$. We now define precisely this catenoid.

\medskip

{\bf The symmetrized catenoid $\Sigma_{ij}$~:} Ideally, we would like to consider for any $\e_{ij} \in (0,1)$, satisfying 
\begin{equation}
|\e_{ij}- \e |\leq \kappa \, \frac{\e}{|\log \e|},
\label{eq:**}
\end{equation}
a {\small CMC}-$1$ surface close to the catenoid parameterized by
\[
X^0 (s, \theta) : = (a^j_i, 1) +  \e_{ij}\, (  \cosh s \, \cos \theta, \cosh s \, \sin \theta, s)
\] 
for $(s, \theta) \in (-\frac{1}{\e_{ij}}, +\infty) \times S^1$, and we would like to ask that this surface is invariant by the inversion centered at $a^j_i$. Using (\ref{eq:2-1}), it is easy to check, as $s$ tends to $+\infty$, the mean curvature of the catenoid parameterized by $X^0$ tends to $1$ (provided the orientation is chosen in such a way that the normal vector points upward when $s$ tends to $+\infty$) while, as $s$ tends to $-\frac{1}{\e_{ij}}$, and as $\e_{ij}$ tend to $0$, its mean curvature tends to $-1$. This is the reason why we will modify the above catenoid by considering  not only the catenoid parameterized by  $X^0$ but also its image by ${\cal I}_{a^j_i}$, the inversion centered at $a^j_i$. This image is parameterized by 
\[
X_0 (s, \theta) := (a_i^j,0) + \frac{1}{\e^2_{ij} \, \cosh^2 s+ (1+\e_{ij} \, s)^2} \,  (\e_{ij}  \, \cosh s \, \cos \theta, \e_{ij} \, \cosh s \, \sin \theta, 1 + \e_{ij} \, s) .
\]
It should be clear that the image of  $[-1, 1]\times S^1$ by $X_0$ can be written as a normal geodesic graph over the surface parameterized by $X^0$ for some function which is bounded by a constant times $\e^2$ in ${\cal C}^{2, \alpha}$ norm and {\it vice-versa}. Hence,  we can easily define a surface of revolution $\Sigma_{ij}$ parameterized by $X_{ij} : {\R}\times S^1 \longrightarrow {\mathbb H}^3$ such that $X_{ij}=X^0$ in  $[1, \infty)\times S^1$, $X_{ij}=X_0$ in $(-\infty, -1]\times S^1$ and the image of $[-1, 1]\times S^1$ by $X_{ij}$ is a geodesic normal graph over the surface parameterized by $X^0$ for some function which is bounded by a constant times $\e^2$ in ${\cal C}^{2, \alpha}$ norm. Moreover,  we can ask that this surface is invariant by ${\cal I}_{a^j_i}$, the inversion centered at $a^j_i$, in the following sense
\[
\forall (s, \theta) \in {\R}\times S^1 , \qquad  \qquad {\cal I}_{a^j_i} \circ X_{ij} (s, \theta) =  X_{ij} (-s, \theta).
\]
In other words ${\cal I}_{a^j_i} \, \Sigma_{ij} =  \Sigma_{ij}$. The advantage of this definition is that,  in order to understand the mean curvature of a surface close to $\Sigma_{ij}$, it is enough to restrict our attention to understanding the mean curvature of a surface close to the upper half of $\Sigma_{ij}$ since we can always reduce to this case by using ${\cal I}_{a^i_j}$. The surface $\Sigma_{ij}$ described above is what we will call a rescaled symmetrized catenoid.

\medskip

{\bf Truncation of the symmetrized catenoid $\tilde \Sigma_{ij}$~: } We shall now explain where we will truncate the symmetrized catenoid.  For all $\rho >0$, we define $s_{ij} > s'_{ij} >0$ by the identities
\[
\e_{ij} \, \cosh s_{ij} =  \rho  \qquad \mbox{and}\qquad   \e_{ij} \, \cosh s'_{ij} = \rho /2.
\]
Observe that, provided $\e$ is small enough,  the set
\[
A^i_j : = {\cal J}^i_j  \, ^{-1}\circ X_{ij} ([s'_{ij}-1, s_{ij}+1]\times S^1),
\]
can be understood as a normal (and hence vertical) geodesic graph over the horosphere ${\cal H}^i(\e)$, for some function $w^i_{j}$, which is defined in some annulus and which is bounded by a constant (depending on $\kappa$) times $\e$ in ${\cal C}^{2, \alpha}$ topology.  In particular, the image of $\del B_\rho (a^j_\infty) $ (resp. $\del B_{\rho/2} (a^j_\infty)$)  first by the vertical geodesic flow and then by ${\cal J}^i_j$ defines on $\Sigma_{ij}$ closed curve $\gamma_j^i $ (resp. $\tilde\gamma_j^i$) which is a constant (depending on $\kappa$) times $\e$ close to $X_{ij} (\{s_{ij}\}\times S^1)$ (resp. $X_{ij} (\{s'_{ij}\}\times S^1)$).  After an inversion with respect to $a_j^i$, a similar analysis can be performed in the region where $\Sigma_{ij}$ is close to ${\cal J}^i_j  ({\cal H}^j(\e))$ and this allows to define a closed curve $\gamma_i^j \subset \Sigma_{ij}$ (resp. $\tilde \gamma ^j_i \subset \Sigma_{ij}$) which is a constant (depending on $\kappa$) times $\e$ close to $X_{ij} (\{- s_{ij}\}\times S^1)$ (resp. $X_{ij} (\{- s'_{ij}\}\times S^1)$).

\medskip

The truncated symmetrized catenoid we will consider is the portion of the symmetrized catenoid which is bounded by the two curves $\gamma_j^i$ and $\gamma^j_i$. This surface will be denoted by $\tilde \Sigma_{ij}$ and its parameterization will be denoted by $\tilde X_{ij} : [- s_{ij}, s_{ij}]\times S^1 \longrightarrow {\HH}^3$.

\medskip

The portion of the truncated symmetrized catenoid between the curves $\gamma^i_j$ and $\tilde \gamma^i_j$ will be parameterized by $(s, \theta) \in [s_{ij}', s_{ij}]\times S^1$ in the following way : first we consider the variables $r : = \e_{ij} \, \cosh s$ and $\theta \in S^1$ as being polar coordinated  in the annulus $ B_\rho (a^j_\infty) - B_{\rho/2} (a^j_\infty)$ and then we take the image of the vertical geodesic graph from the point of coordinates $(s, \theta)$ in $ B_\rho (a^j_\infty) - B_{\rho/2} (a^j_\infty)$ by ${\cal J}^i_j$ to obtain a point in $\Sigma_{ij}$.  The image of $[1 - s_{ij}', s_{ij}' -1]\times S^1$ by $X_{ij}$ will be parameterized as before by $X_{ij}$, and finally we will use a parameterization which interpolates between the two above defined parameterization in the intermediate region. 

\medskip

{\bf A transverse vector field $\tilde N_{ij}$~:} We end this section by defining, at every point $ p \in \tilde X_{ij} ([s_{ij}-1, s_{ij}]\times S^1)$, a transverse vector field $p \rightarrow \tilde N (p)$ which is the tangent vector at $p$ of the unique normal geodesic issued from some point of ${\cal J}^i_j \, {\cal H}^i (\e)$, which meets $\tilde \Sigma_{ij}$ at the point $p$. On $\Sigma_{ij}$, the Euclidean normal vector is chosen to be the inward pointing normal vector, hence, it is given by
\[
N_{eucl} (s, \theta) : = \frac{1}{\cosh s}\, (- \cos \theta, - \sin \theta,  \sinh s),
\]
for all $s > 1$. Now,  we define $\tilde N_{ij}$ a vector field  on $\tilde \Sigma_{ij}$ which is everywhere transverse to $\tilde \Sigma_{ij}$ such that $\tilde N_{ij} = N_{eucl}$ for all $s \in [1, s_{ij}' -1]$ and $\tilde N_{ij}$ coincides with $\tilde N$  for all $s \in [s'_{ij}, s_{ij}]$. 

\section{Constant mean curvature $1$ surfaces close a horosphere}

In this section, we find {\small CMC}-$1$ surfaces with boundary, which are complete, embedded and which are close to a horosphere with finitely many balls removed. To this aim, we rephrase the problem in term of a nonlinear partial differential equation problem which we solve by using the implicit function theorem in a suitably defined function space. So to begin with we briefly develop the relevant linear analysis, then turn to the study of the nonlinear operators which appear in the expression of the mean curvature operator for a graph. 

\subsection{Function spaces and linear analysis}

To begin with, let us define the function space we will work with in this section. 
\begin{definition}
Given $k \in {\N}$, $\alpha \in (0,1)$, $\delta \in {\R}$, we define the weighted H{\"o}lder space ${\cal C}^{k, \alpha}_\delta ({\R}^2)$ to be the space of functions in ${\cal C}_{loc}^{k, \alpha} ({\R}^2)$  for which the following norm is finite
\[
\|u\|_{{\cal C}^{k, \alpha}_\delta} := [u]_{k, \alpha, \overline B_1} +  \sup_{r >1/2}  r^{-\delta} \,  [u (r \,  \cdot)]_{k, \alpha, B_2 -B_1} ,
\]
where $[u]_{k, \alpha , \Omega}$ is the usual ${\cal C}^{k, \alpha}$ H{\"o}lder norm in $\Omega$.
\label{de:4}
\end{definition}
For any bounded  open set $\Omega \subset {\R}^2$, we denote  by ${\cal C}^{k,\alpha}_\delta({\R}^2-\Omega)$ the restriction  of ${\cal C}^{k,\alpha}_\delta({\R}^2) $ to ${\R}^2-\Omega$. This space is naturally endowed with the induced  norm. We will denote by  by $[{\cal C}^{k,\alpha}_{\delta}({\R}^2-\Omega)]_0$ the subspace of ${\cal C}^{k,\alpha}_\delta({\R}^2-\Omega)$ corresponding to functions vanishing on $\del \Omega$.

\medskip

First, observe that the operator 
\[
\Delta : [{\cal C}^{2, \alpha}_{\delta } ({\R}^2 - \Omega)]_0 \longrightarrow {\cal C}^{2, \alpha}_{\delta -2} ({\R}^2 - \Omega) ,
\]
is well defined and bounded. We shall denote this operator by $\Delta_\delta$, the subscript $\delta$ referring to the weight, and hence to the function space where the Laplacian is defined. It follows from the general theory of  elliptic partial differential operators with regular singularities that the operator $\Delta_\delta$  is Fredholm for all $\delta \notin {\mathbb Z}$. It is also well known that, if $\delta \in {\R}-{\mathbb Z}$, then $\Delta_\delta$ is injective if and only if  $\Delta_{-\delta}$ is surjective. Furthermore, the dimension of the kernel of $\Delta_{-\delta}$ is equal to the dimension of the cokernel of $\Delta_\delta$. In our case, this reduces to the following result which classifies the range of weight parameters $\delta$ for which this operator is injective or surjective.  
\begin{proposition}
For any $m \in {\mathbb N}$. Assume that $-m-1 < \delta <-m$, then $\Delta_\delta$ is injective but not surjective and in fact has a cokernel of dimension $m+1$.  
 Assume that $m < \delta < m+1$, then $L_\delta$ is surjective but not injective and in fact has a kernel of dimension $m + 1$.  
\end{proposition}
When the weight parameter is negative, we can still make the operator $L_\delta$ surjective by considering a finite dimensional extension of the above defined H{\"o}lder spaces.  To this aim, let us fix $R_1 >0$ such that $\Omega \subset B_{R_1}$.  We will need a radial cutoff function $\chi$ which is identically equal to $1$ away from $B_{4R_1}$  and which is identically equal to $0$ in $B_{2R_1}$. We define
\[
{\cal D}_0 : = \mbox{Span}  \{ \chi  \,,  \chi \, \log r\}.
\]
This space will be identified with ${\R}^2$ and will be  induced with the Euclidean norm.

\medskip

For all $\delta \in (-1,0)$, we claim that, not only $\mbox{Ker}\,  \Delta_{-\delta} \subset  [{\cal C}^{2, \alpha}_{-\delta } ({\R}^2 - \Omega)]_0$ but we even have
\[
\mbox{Ker}\,  \Delta_{-\delta} \subset  [{\cal C}^{2, \alpha}_{\delta } ({\R}^2 - \Omega)]_0 \oplus {\cal D}_0 .
\]
Indeed, we can expend any solution of $\Delta w=0$ as
\[
w = \sum_{m \in {\mathbb Z}} w_m  \, e^{im\theta}. 
\]
It is easy to see that, if $w \in {\cal C}^{2, \alpha}_{-\delta } ({\R}^2 - \Omega)$, then 
\[
\forall m \neq 0, \qquad \qquad |w_m | \leq c\, r^{-|m|},
\]
in ${\R}^2 - B_{R_1}$, for some constant $c$ which does not depend on $m \neq 0$. Moreover $w_0$ is a linear combination of $1$ and $\log r$. The claim then follows at once.

\medskip

Using this, we can prove a "Linear Decomposition Lemma" as in \cite{Kus-Maz-Pol} or in \cite{Maz-Pol-Uhl} (see also \cite{Pac-Riv} for more details on this kind of decomposition)~:
\begin{proposition} 
Assume that $\delta \in (-1,0)$. We can decompose
\[
{\cal D}_0 = {\cal N}_0 \oplus {\cal K}_0 ,
\]
where ${\cal K}_0 \neq \mbox{Span}\{\chi\}$ and ${\cal N}_0$ are both of dimension $1$, in such a way that
\[
\mbox{Ker} \, \Delta_{-\delta} \subset  [{\cal C}^{2, \alpha}_{\delta } ({\R}^2 - \Omega)]_0 \oplus {\cal N}_0 .
\]
and the operator 
\[
\Delta : [{\cal C}^{2, \alpha}_{\delta } ({\R}^2 - \Omega)]_0 \oplus {\cal K}_0  \longrightarrow  {\cal C}^{0, \alpha}_{\delta -2} ({\R}^2- \Omega ) , 
\]
is an isomorphism. 
\label{pr:4}
\end{proposition}
{\bf Proof~:} To begin with, we claim that
\[
\Delta : [{\cal C}^{2, \alpha}_{\delta } ({\R}^2 - \Omega)]_0 \oplus \mbox{Span} \{\chi\}  \longrightarrow  {\cal C}^{0, \alpha}_{\delta -2} ({\R}^2- \Omega ) , 
\]
is an isomorphism. Let us denote by $G$ its inverse. Granted this, we see that the kernel of  
\[
\Delta : [{\cal C}^{2, \alpha}_{\delta } ({\R}^2 - \Omega)]_0 \oplus {\cal D}_0  \longrightarrow  {\cal C}^{0, \alpha}_{\delta -2} ({\R}^2- \Omega ) , 
\]
is one dimensional and, in fact, is spanned by the function $x \longrightarrow \chi \, \log r - G \, (\Delta \, (\chi \, \log r) )$. This shows that we can find ${\cal N}_0$, a one dimensional subspace of ${\cal D}_0$, in such a way that 
\[
\mbox{Ker} \, \Delta_{-\delta} \subset  [{\cal C}^{2, \alpha}_{\delta } ({\R}^2 - \Omega)]_0 \oplus {\cal N}_0 .
\]
Finally, we can choose any one dimensional subspace ${\cal K}_0 \neq \mbox{Span}\{\chi\}$ such that ${\cal D}_0 = {\cal N}_0 \oplus {\cal K}_0$ and, using $G$, one can prove easily that
\[
\Delta : [{\cal C}^{2, \alpha}_{\delta } ({\R}^2 - \Omega)]_0 \oplus {\cal K}_0  \longrightarrow  {\cal C}^{0, \alpha}_{\delta -2} ({\R}^2- \Omega ) , 
\]
is an isomorphism.

\medskip

It remains to prove the claim. We first solve the equation $\Delta v =f$ in ${\R}^2 - B_{R_1}$, with boundary condition $v \in {\R}$ on $\del B_{R_1}$. To this aim we consider the Fourier series of both $f$ and $v$
\[
f = \sum_{m \in {\mathbb Z}}  f_m \, e^{i m\theta} \qquad \mbox \qquad v = \sum_{m \in {\mathbb Z}}  v_m \, e^{i m\theta}.
\]
For each $m \in {\mathbb Z}$, the function $v_m$ has to solve the ordinary differential equation
\begin{equation}
\del^2_r v_m + \frac{1}{r}\, \del_r v_m - \frac{n^2}{r^2}\, v_m =f_m .
\label{eq:NBNB}
\end{equation}
Since $|f_m(r) |\leq \|f\|_{{\cal C}^{0, \alpha}_{\delta -2}}\, r^{\delta -2}$ in ${\R}^2 - B_{R_1}$, we see that, for all $m \neq 0$, the function 
\[
r \longrightarrow \frac{\|f\|_{{\cal C}^{0, \alpha}_{\delta -2}}}{m^2 - \delta} \, r^{\delta -2},
\]
can be used as a barrier function and the method of sub and supersolutions allows one to prove the existence of $v_m$, solution of (\ref{eq:NBNB}) with $v_m (R_1)=0$, which satisfies
\[
|v_m(r) |\leq \frac{\|f\|_{{\cal C}^{0, \alpha}_{\delta -2}}}{m^2- \delta^2} \, r^{\delta}.
\]
Finally, for $m=0$, we simply define
\[
v_0(r)  =  \int_r^\infty t^{-1}\, \int_{t}^\infty s \, f_0 (s) \, ds \, dt .
\]
It is then a simple exercise to show that
\[
\sup_{{\R}^2 - B_{R_1}} r^{-\delta} \, |v |\leq c \, \|f\|_{{\cal C}^{0, \alpha}_{\delta-2}} ,
\]
which, together with Schauder's estimates, yields
\[
\| \chi \, v \|_{{\cal C}^{2, \alpha}_\delta} \leq c \, \|f\|_{{\cal C}^{0, \alpha}_{\delta-2}}.
\]
Observe that $\tilde f  : = f - \Delta (\chi \, v)$ is supported in $B_{4R_1} - \Omega$. 

\medskip

In order to solve $\Delta \tilde v =  \tilde f $ in ${\R}^2 - \Omega$, we first perform a Kelvin transform. That is, we define
\[
\bar v (x) : = \tilde v \left( \frac{x}{|x|^2}\right) \qquad \mbox{and}\qquad  \bar f (x) : = \frac{1}{|x|^4} \, \tilde f \left( \frac{x}{|x|^2}\right),
\]
and reduce the problem to solve $\Delta \bar v = \bar f$ in some bounded set $\bar \Omega$, with boundary data $\bar v = 0$. It is well known that $\bar v$ exists and that we have the estimates
\[
|\bar v (x) - \bar v(0)|\leq c\, r \, \| f \|_{{\cal C}^{0, \alpha}_{\delta -2}}\qquad \mbox{and}\qquad   |\bar v(0)|\leq c \, \| f \|_{{\cal C}^{0, \alpha}_{\delta -2}}.
\]
Performing the Kelvin transform backward, and making use of Schauder's estimates, we conclude that
\[
\| \tilde v - \bar v(0) \|_{{\cal C}^{2, \alpha}_{-1}} \leq c \, \|f\|_{{\cal C}^{0, \alpha}_{\delta-2}}.
\]
Observe that  $v + \tilde v$ is a solution of our problem which belongs to $[{\cal C}^{2, \alpha}_{\delta } ({\R}^2 - \Omega)]_0 \oplus \mbox{Span} \{\chi\}$. The ends the proof of the claim.   \hfill $\Box$

\begin{remark}
The choice of ${\cal K}_0 \neq \mbox{Span}\{\chi\}$ implies that, for all $\delta \in (-1, 0)$,  the operator 
\[
\Delta : {\cal C}^{2, \alpha}_{\delta } ({\R}^2) \oplus {\cal K}_0  \longrightarrow  {\cal C}^{0, \alpha}_{\delta -2} ({\R}^2) , 
\]
is injective. This property will be used in \S 6.
\label{re:inj}
\end{remark}

The last result we would like to mention is concerned with the asymptotic behavior of solutions of the homogeneous problem $\Delta w =0$ which belong to the space ${\cal C}^{2, \alpha}_\delta ({\R}^2 -\Omega)\oplus {\cal D}_0$, an easy proof can be obtained by considering once more  the Fourier decomposition of $w$.  
\begin{proposition}
Let $\delta \in (-1, 0)$ be fixed. There exists a constant $c >0$ such that, all $w\in {\cal C}^{2, \alpha}_\delta ({\R}^2 -\Omega) \oplus {\cal D}_0$, solution of $\Delta w=0$ in ${\R}^2 -\Omega$,  can be decomposed as $w =  v+  a \, \chi \, \log r + b \, \chi$ with 
\[
|a|+ |b|+ \| v \|_{{\cal C}^{2, \alpha}_{ -1}}\leq c \,  \|w\|_{{\cal C}^{2, \alpha}_\delta \oplus {\cal D}_0}.
\]
\label{pr:5}
\end{proposition}
        
\subsection{Some nonlinear differential operators}

We introduce some second order nonlinear differential operators and study their properties. These operators appear in the expression of the mean curvature operator of a vertical graph in the upper half space model of ${\HH}^3$. To begin with, we set
\[
K(u) : = u \, \Delta u - |\nabla u|^2.
\]
Observe that, for any $t\in {\R}$, we have
\[
K(r^t \, u) = r^{2t}\, K(u) .
\]
We also have the following easy~:
\begin{lemma}
For all, $\delta \leq 0$, the nonlinear mapping 
\[
{\cal C}^{2, \alpha}_\delta ({\R}^2 - \Omega) \ni u  \longrightarrow  K(u) \in {\cal C}^{0, \alpha}_{\delta-2}({\R}^2 -\Omega) ,
\]
is ${\cal C}^\infty$.
\label{le:2}
\end{lemma}
Now, we define the first order partial differential operator
\[
Q_1 (u) :=  2 \, \left( 1 + \frac{3}{2} \, |\nabla u|^2 - (1 + |\nabla u|^2)^{\frac{3}{2}}\right) .
\]
Using the fact that the function $\zeta \longrightarrow 1 + \frac{3}{2} \, \zeta^2  - (1+\zeta^2)^{\frac{3}{2}}$ is ${\cal C}^\infty$ and has all its derivatives up to order $3$ which vanish at $0$, we get 
\begin{lemma}
For all $t_0 \in {\R}$ and for all $\delta > 2 \, t_0 -2$, the nonlinear operator
\[
(-\infty, t_0) \times  {\cal C}^{2, \alpha}_0 ({\R}^2 -B_{R_1}) \ni (t, u) \longrightarrow r^{-2t} \, Q_1 (r^t \, (1+u) ) \in{\cal C}^{0, \alpha}_{\delta-2}({\R}^2 -B_{R_1}),
\]
is ${\cal C}^\infty$.
\label{le:3}
\end{lemma}
Finally, we introduce  the second order nonlinear partial differential operator
\[
Q_2 (u)  : = u \, \Delta u \, |\nabla u|^2  -  \frac{u}{2} \,  \nabla u \cdot \nabla |\nabla u|^2  
\]
This expression is quartic in $u$ and its derivatives and it is easy to check that 
\begin{lemma}
For all $t_0\in {\R}$ and for all $\delta > 2 \, t_0 -2$, the nonlinear operator
\[
(-\infty , t_0) \times  {\cal C}^{2, \alpha}_0  ({\R}^2 - B_{R_1}) \ni (t, u) \longrightarrow r^{-2t} \, Q_2 (r^t \,(1+u) ) \in{\cal C}^{0, \alpha}_{\delta-2}({\R}^2 -B_{R_1}) ,
\]
is ${\cal C}^\infty$.
\label{le:4}
\end{lemma}

\subsection{Application of the implicit function theorem to perturb ${\cal H}^i$}

We now fix $\delta \in (-1, 0)$. We choose $\rho >0$ small enough so that, for $j \in I_i$, the $B_{2 \, \rho}(a^j_\infty)$ are disjoint. We now define 
\[
\Omega_{i} : =  \cup_{j=1}^{n_i} B_{\rho} (a^j_\infty),
\]
where $n_i$ is the cardinal of $I_i$. Given functions $\phi^i_j \in {\cal C}^{2, \alpha} (S^1)$,  we set $\Phi_i : = (\phi^i_j)_{j \in I_i}$ and define $W_{\Phi_i}$ to be the unique harmonic extension of the functions 
\[
 \del B_\rho (a^j_\infty)\ni x  \longrightarrow \phi^i_j \left( \frac{x- a^j_\infty}{\rho}\right),
\]
which belongs to ${\cal C}^{2, \alpha}_\delta({\R}^2-\Omega_{i}) \oplus {\cal K}_0$. This solution is obtained using the result of Proposition~\ref{pr:4} and thanks to Proposition~\ref{pr:5}, we know that $W_{\Phi_i}$ can be decomposed as
\[
W_{\Phi_i} = v_{\Phi_i} +  a_{\Phi_i} \, \chi \, \log r + b_{\Phi_i} \, \chi 
\]
with
\[
|a_{\Phi_i}|+ |b_{\Phi_i}| + \| v_{\Phi_i} \|_{{\cal C}^{2, \alpha}_{-1}}\leq c \, \|\Phi_i \|_{{\cal C}^{2, \alpha}}.
\]
Any function $w \in [{\cal C}^{2, \alpha}_\delta ({\R}^2 -\Omega_{i})]_0 \oplus {\cal K}_0$, can be uniquely decomposed as 
\[
w = v +  a \, \chi \, \log r + b \, \chi  ,
\]
with obvious notations. With this decomposition in mind, we can define 
\[
U (w , {\Phi_i}) : =   \chi \, r^{a_{\Phi_i} +a }\,  (1+b_{\Phi_i} + b + v_{\Phi_i} + v) + (1-\chi) \, (1+ v_{\Phi_i} +  v)
\]
where $\chi$ is the cutoff function which has been used in \S 4.1 to define the space ${\cal D}_0$. 

\medskip

We define ${\cal E}$ to be the set of functions 
\[
(w, \Phi_i) \in  \left( [{\cal C}^{2, \alpha}_{\delta} ({\R}^2 -\Omega_{i})]_0  \oplus {\cal K}_0 \right) \times [{\cal C}^{2, \alpha}(S^1)]^{n_i} ,
\]
such that the function $w + W_{\Phi_i}$ can be decomposed as
\[
w+ W_{\Phi_i} = v +  a \, \chi \, \log r + b \, \chi 
\]
with $v \in {\cal C}^{2, \alpha}_\delta ({\R}^2 -\Omega_{i})$, $b  \in {\R}$ and $a < \frac{2+\delta}{2}$. Recall that $n_i$ is the cardinal of $I_i$. We also define 
\[
{\cal F} : = {\cal C}^{0, \alpha}_{\delta-2} ({\R}^2-\Omega_{i}) .
\] 
It follows directly from Lemma~\ref{le:2}- Lemma~\ref{le:4} that the mapping
\[
{\cal N} :  (w, \Phi_i )  \in {\cal E}  \longrightarrow  (1+r^2)^{-(a_{\Phi_i} +a)} \,  \left( K(U)+  Q_1 (U ) + Q_2 (U) \right) \in {\cal F},
\]
is well defined and of class ${\cal C}^\infty$. Furthermore, $D_w {\cal N}$, the differential of ${\cal N}$ with respect to $w$, computed at $w=0$ and  $\Phi_i = 0$, is given by $\Delta$ and, thanks to Proposition~\ref{pr:4}, it is an isomorphism. We can apply the implicit function theorem to obtain,  for all $\Phi_i $ small enough,  the existence of $w_{\Phi_i}$ solution of ${\cal N} (w, {\Phi_i})=0$. 

\medskip

Observe that the equation $H(u)=1$, where $H(u)$ is defined in (\ref{eq:2-2}) reduces to 
\[
K(u)+ Q_1(u)+Q_2(u)=0.
\]
Hence, the graph of $U(w_{\Phi_i}, \Phi_i)$ produces a {\small CMC}-$1$ surface, denoted by $M_i (\Phi_i)$, whose boundary is parameterized by $\Phi_i$ and which has a regular end asymptotic to the end of a catenoid cousin. It will be important to observe that the following expansion holds
\[
U (w_{\Phi_i}, \Phi_i)  = 1 - \e \, \log \e + W_{\Phi_i} + V_{\Phi_i} ,
\] 
in any $B_{2\rho}(a^j_\infty)- B_\rho (a^j_\infty)$. Furthermore, if we assume that
\[
\|\Phi_i\|_{{\cal C}^{2, \alpha}}\leq \kappa \, \e \qquad \mbox{and}\qquad  \| \tilde \Phi_i\|_{{\cal C}^{2, \alpha}}\leq \kappa \, \e , 
\]
we have the estimate
\begin{equation}
\| V_{\Phi_i} - V_{\tilde \Phi_i}\|_{{\cal C}^{2, \alpha}(B_{2\rho}-B_\rho)} \leq c \, \e \, \| \Phi_i - \tilde \Phi_i\|_{{\cal C}^{2, \alpha}} ,
\label{eq:cont-1}
\end{equation}
for some constant $c >0$ only depending on $\kappa$.

\section{Mean curvature $1$ surfaces close to a vertical rescaled symmetrized catenoid}

In this section we would like to find {\small CMC}-$1$ surfaces close to a truncated,  rescaled, symmetrized vertical catenoid. The strategy is very close to what we have done in the previous section. To begin with, we define the function spaces and develop the necessary linear analysis. Then we rephrase our problem into a  the nonlinear partial differential equation which, this time,  we solve by some contraction mapping argument.

\subsection{Function spaces and linear analysis}

To begin with, let us define the function space we will work with in this section. 
\begin{definition} 

Given $k \in {\N}$, $\alpha \in (0,1)$ and $\mu  \in {\R}$, we define
the 
weighted Holder space ${\cal C}^{k, \alpha}_\mu ({\R}\times S^1)$ to be the space of functions in ${\cal C}_{loc}^{k, \alpha} ({\R}\times S^1)$  for which the following norm is finite
\[
 \ ||w||_{{\cal C}^{k,\alpha}_\mu}: = \sup_{s \in {\R}} \,  (\cosh s)^{-\mu}\,   [w]_{k,\alpha,[s,s+1]\times S^1} 
\]
where $[w]_{k,\alpha, [s, s +1]\times S^1}$ denotes the usual ${\cal C}^{k,\alpha}$ H{\"o}lder norm on the  set $[s,s+1] \times S^1$. 
\label{de:5}
\end{definition} 

For any closed interval $I\subset\R$, we denote the restriction  of ${\cal C}^{k,\alpha}_\mu({\R} \times S^1) $ to $I\times S^1$ by ${\cal C}^{k,\alpha}_\mu(I \times S^1)$, endowed with the induced  norm.

\medskip

Let us introduce the operator
\begin{equation}
{\cal L} : = \del^2_s + \del_\theta^2 + \frac{2}{\cosh^2 s} .
\label{eq:5-1}
\end{equation}
The mapping properties of ${\cal L}$ which will be needed are contained in the~:

\begin{proposition}
Fix $\mu \in (1,2)$. Then for any $s_0 >1$ there  exists an operator
\[
{\cal G}_{s_0} : {\cal C}^{0,\alpha}_{\mu}([-s_0, s_0]\times S^1) \longrightarrow  {\mathcal C}^{2,\alpha}_{\mu}([-s_0, s_0] \times S^1) ,
\]
such that for any $f \in {\cal C}^{0,\alpha}_{\mu}  ( [-s_0, s_0] \times S^1)$, the function $w={\cal G}_{s_0}(f)$ solves $ {\cal L} w = \ds  f $ in $(- s_0,s_0) \times S^1$ with $w \in \mbox{Span} \{1, e^{i \theta}, e^{-i\theta}\}$ on $\{\pm s_0\} \times  S^1$. Moreover, $||{\cal G}_{s_0}(f)||_{{\cal C}^{2, \alpha}_\mu} \leq c \, ||f||_{{\cal C}^{0, \alpha}_\mu },$ for some constant $c>0$ independent of $s_0$.
\label{pr:6}
\end{proposition}
Before, we proceed to the proof of this result, let us briefly comment on its statement.  Observe that we have restricted $s_0$ to be larger than $1$. This condition is needed in order to be able to apply Schauder's estimate with some constant which is independent of $s_0$. Moreover, we have not asked that the boundary data should be equal to $0$ but rather should belong to the space spanned by the first three eigenfunctions of the Laplacian on $S^1$. This is needed to guaranty that we can find a right inverse for ${\cal L}$ whose norm is independent of $s_0$. Finally, observe that we do not have uniqueness of ${\cal G}_{s_0}$. However, in the forthcoming analysis we will always use the right inverse constructed in this Proposition.

\medskip

\noindent
{\bf Proof:} We decompose both $w$ and $f$ into Fourier series
\[
w= \sum_{m \in {\Z}} w_m \, e^{im \theta} \qquad \mbox{and} \qquad  f= \sum_{m \in {\Z}} f_m \,  e^{im \theta}.
\]
Then $w_m$ must solve
\[
\del^2_s w_m -m^2 w_m + \frac{2}{\cosh^2 s} w_m = f_m  .
\]
For  $|m|\geq 2$,
\[
L_m = \del_s^2 - m^2  + \frac{2}{\cosh^2 s}
\]
satisfies the maximum principle, so that if $w$ is defined on some interval  $[s_1,s_2]\subset {\R}$ and if $w(s_1) \geq 0$, $w(s_2) \geq 0$ and  $L_m w \leq 0$ on  $(s_1,s_2)$, then $w \geq 0$ in $[s_1,s_2]$.  We obtain the solution of $L_m w_m = f_m$ by the method of sub and supersolutions once we have constructed an appropriate barrier  function. But
\[
\begin{array}{rlll}
L_m (\cosh s)^{\mu} & = & \ds \left(  (\mu^2 - m^2) \cosh^2s + 2+\mu -  \mu^2 \right) (\cosh s)^{\mu-2}\\[3mm]
                     &  \leq &  - \, (m^2 -2 - \mu) \,  (\cosh s)^{\mu} .
\end{array}
\]
Hence, using the method of sub and super solutions, we have both the existence of $w_m$ and the estimate
\[
\sup_{[-s_0, s_0]}   (\cosh s)^{-\mu} \, |w_m | \leq (m^2-2-\mu)^{-1} \sup_{[-s_0, s_0]} (\cosh s)^{-\mu} \, |f_m|.
\]
Next we obtain the solution and estimates when $m=0, \pm 1$.  These
solutions are explicitely given by
\[
w_{0}(s) = \tanh s \int_0^s \tanh^{-2} t \int_0^t \tanh u \, f_{0}(\xi)\, d\xi \, dt,
\]
and
\[
w_{\pm 1} (s) = \cosh^{-1} s \int_0^s \cosh^2 t \int_0^t \cosh^{-1} \xi \, f_{\pm 1}(\xi)\, d\xi\, dt.
\]
Straightforward estimates using these formul\ae\ yield
\[
\sup_{[-s_0, s_0]}  (\cosh s)^{-\mu} \,  (|w_0 | + |w_{\pm 1}|) \leq c \,  \sup_{[-s_0, s_0]}  (\cosh s)^{-\mu} \, (|f_0 | + |f_{\pm 1}|) ,
\]
for some constant $c>0$ independent of $s_0$. Summing over $m$ we conclude that 
\[
\sup_{[-s_0, s_0]\times S^1} (\cosh s)^{-\mu} \, |w|\leq c \, \sup_{[-s_0, s_0]\times S^1} (\cosh s)^{-\mu} \, |f|
\]
The estimates for the  derivatives of $w$ are then obtained by Schauder theory.  \hfill $\Box$

\medskip

We will also need some properties of the Poisson operator for the Laplacian on a half cylinder. More precisely, we define
$$
{\cal P}:  {\cal C}^{2,\alpha}(S^1) \longrightarrow {\mathcal C}^{2,\alpha}_{2}([0, \infty)  \times S^1) ,
$$
such that for all $\phi  \in {\cal C}^{2,\alpha}(S^1)$ the function $w : = {\cal P}(\phi)$ is bounded,  solves $(\del_s^2 + \del_\theta^2) \, w =0$ in $(0, \infty)\times S^1$ and satisfies $w=\phi$ on $\{0\}\times S^1$.
\begin{proposition}
There exists $c>0$ such that for all $\phi  \in {\cal C}^{2,\alpha}(S^1)$ which is orthogonal to $1$, $e^{i\theta}$ and $e^{-i\theta}$  in the $L^2$ sense on $S^1$, we have  $||{\cal P}(\phi)||_{2,\alpha,-2} \leq c \, ||\phi ||_{2,\alpha}$.
\label{pr:7}
\end{proposition}
{\bf Proof:} Again, we decompose $\phi$ into Fourier  series
$$
\phi = \sum_{|m|\geq 2} \phi_m  \, e^{im \theta} ,
$$
and solve explicitely
$$
w= \sum_{|m|\geq 2} \phi_n \, e^{-|m|s} \, e^{im \theta} .
$$
From this it is easy to get the estimate $\sup_{[1, \infty)\times S^1} (\cosh s)^{2} \, |w | \leq c \, \|\phi\|_{L^\infty}$. Now, it suffices to apply the maximum principle to estimate $w$ in $[0,1]\times S^1$. Estimates for higher order derivatives follow from Schauder's estimates. \hfill $\Box$

\subsection{Structure of the mean curvature operator for surfaces close to a vertical catenoid}

To begin with, we fix $\e_{ij} \in (0,1)$ satisfying (\ref{eq:**}) and consider the truncated symmetrized catenoid $\tilde \Sigma_{ij}$ parameterized by $\tilde X_{ij}$ as described in section \S 3. Recall that we have also defined a transverse vector field $\tilde N_{ij}$. We shall now explain how we will parameterize a surface close to $\tilde \Sigma_{ij}$. This parameterization has to be coherent with our desire to keep everything invariant under the inversion ${\cal I}_{a^j_i}$.  Indeed, we will describe any surface close enough to the image of $[0, s_{ij}] \times S^1$ by $\tilde X_{ij}$ as a geodesic graph over $\tilde \Sigma_{ij}$ of a (small) function $w$, using the vector field $\tilde N_{ij}$. Namely
\begin{equation}
Z_w \, (s, \theta) := \Gamma ( \tilde X_{ij} (s, \theta), \tilde N_{ij} (s, \theta), w  (s, \theta))  \in ({\mathbb H}^3, g_{hyp}).
\label{eq:kk}
\end{equation}
In order to describe surfaces close to the image of $[-s_{ij}, 0]\times S^1$ by $\tilde X_{ij}$, we first perform an inversion with respect to $a^j_i$ and reduce the problem to the former case.  We are interested in computing the mean curvature of the surface parameterized by $Z_w$. We will not need the exact expression for the mean curvature operator in terms of the function $w$ but rather we need to understand its structure. This is the content of the following~:
\begin{proposition}
The surface parameterized by $Z_w$  has mean curvature $1$ if and only if the function $w$ is a solution of 
\begin{equation}
{\cal L} w  =  f_{ij}  + \e_{ij} \, L w  + (\e_{ij}+ \e_{ij}^2 \, \cosh^2 s)  \, Q_2 \left(\frac{w}{\e_{ij} \, \cosh s}\right)+ \e_{ij} \, \cosh s \, Q_3 \left( \frac{w}{\e_{ij} \, \cosh s }\right)  ,
\label{eq:nn}
\end{equation}
where 
\[
\|f_{ij}\|_{0, \alpha, 0}\leq c \, \e_{ij}^2 ,
\]
and in fact 
\[
f_{ij} = \e_{ij}^2 \, \cosh^2 s \, \frac{1- \tanh s}{1+\e_{ij} \, s} ,
\]
when $1 \leq s \leq s'_{ij} -1$. Where ${\cal L}$ is the linear second order differential operator which has already been defined in (\ref{eq:5-1}), $L$ is a linear second order differential operator,  $Q_2$,  $Q_3$ are nonlinear second order differential operators which satisfy
\[
Q_2 (0) =  Q_3 (0) = 0 \qquad   D_w Q_2 (0) = D_w Q_3 (0) =0 \qquad \mbox{and} \qquad  D^2_w Q_3 (0) =0. 
\]
Moreover the coefficients of $L$ and the coefficients of the Taylor expansion of $Q_2$ and $Q_3$ are functions of $s$  which are bounded uniformly in $s$, as are all of their derivatives, independently of $\e_{ij} \in (0,1)$, though they depend on $\kappa$.
\label{pr:8}
\end{proposition}
{\bf Proof :} The proof of this proposition is close to the proof of the similar result in Euclidean space which is given in \cite{Maz-Pac-Pol-1}. However, since this technical result is a key point of our construction we briefly sketch the proof. For the sake of simplicity, we omit the $i,j$ indices and we set
\[
\tilde{w} : = \ds \frac{w}{\e \, \cosh s}.
\]

\medskip

We will restrict our attention to the case where $s \in [1, s'_{ij}-1]$. Obvious modifications are needed to treat the case where $s \in [0,1]$ and $s\in [s'_{ij}, s_{ij}]$. In order to compute the mean curvature of the surface parameterized by $Z_w$ we first consider this surface to be embedded in ${\R}^3$ and compute the coefficients of the first and second fundamental form. A simple computation shows that the coefficients of the first fundamental  form of the surface parameterized by (\ref{eq:kk}) are given by 
\[
E(w): = \del_s Z_w \cdot \del_s Z_w  = \e^2 \, \left( \cosh^2 s - 2 \cosh s \, \tilde{w}  +  \cosh^2 s \, P_E (\tilde w , \del_s \tilde w) \right), 
\]
\[
F(w) : = \del_s Z_w \cdot \del_\theta Z_w  =  \e^2 \, P_F (\tilde w , \del_s \tilde w, \del_\theta \tilde w)  ,
\]
and
\[
G(w) : = \del_\theta Z_w \cdot \del_\theta Z_w =   \e^2 \, \left( \cosh^2 s + 2 \cosh s \, \tilde{w}  +   \cosh^2 s \, P_G (\tilde{w}, \del_\theta \tilde w )\right). 
\]
Here $P_E, P_F$ and $P_G$ fulfill properties similar to the those enjoyed by $Q_2$ in the statement of the result.  Collecting this, we obtain
\[
(E G -F^2) (w) = \e^4 \, \cosh^4 s \,  \left( 1  + P_{EG-F^2}  (\tilde{w}, \nabla \tilde{w}) \right),
\]
where $P_{EG-F^2}$ fulfills properties similar to the those enjoyed by $Q_2$ in the statement of the result. 

\medskip

In the same way, we compute the coefficients of the second fundamental form and find that these are given by 
\[
\begin{array}{rllll}
\sqrt{E G -F^2} \, e (w)  & :  = & \del_s Z_w \times \del_\theta Z_w \cdot \del_s^2 Z_w \\[3mm]
                          & = &  -\e^3 \, \cosh^2 s \, \left(1 + \cosh s \,  \del_s^2 \, \tilde{w} + \sinh s \, \del_s \tilde{w} + P_e(\tilde{w},\nabla\tilde{w}, \nabla^2 \tilde{w})\right),
\end{array}
\]
\[
\begin{array}{rllll}
\sqrt{EG -F^2} \,  f (w) & : = & \del_s Z_w \times \del_\theta Z_w \cdot \del_s \del_\theta Z_w \\[3mm]
                        & =   & - \e^3 \, \cosh^2 s  \,  \left( \cosh s \, \del_s \del_\theta \tilde{w} +  \, P_f (\tilde{w}, \nabla \tilde{w},  \nabla^2 \tilde{w})  \right) ,
\end{array}
\]
and 
\[
\begin{array}{rlllll}
\sqrt{E G -F^2} \,  g (w) & : = & \del_s Z_w \times \del_\theta Z_w \cdot \del_\theta^2 Z_w \\[3mm]
                         &   = &  \ds - \e^3 \, \cosh^2 s \, \left(- 1 + \cosh s \,  \del^2_\theta \tilde{w} + \left( \cosh s - \frac{2}{\cosh s} \right) \, \tilde{w}  \right. \\[3mm]
                         &       & \left. + \sinh s \, \del_s \tilde{w} + P_g (\tilde{w}, \nabla \tilde{w},  \nabla^2 \tilde{w}) \right) ,
\end{array}
\]
where, here also, $P_e, P_f, P_g$ fulfill properties similar to the those enjoyed by $Q_2$ in the statement of the result. 

\medskip

The Euclidean mean curvature operator may then be expressed in terms of these coefficients as 
\[
H_{eucl}(w)  : = \frac{1}{2} \ \frac{e_w G_w - 2 f_w F_w + g_w E_w}{E_w G_w -F_w^2}.
\]
Using the previous expansions we obtain
\[
H_{eucl} (w) = \frac{1}{\e^2\, \cosh^2 s} \, {\cal L} \, w + \frac{1}{\e \, \cosh^2 s} \, \tilde Q_2 \left( \frac{w}{\e \, \cosh s}\right) + \frac{1}{\e \, \cosh s}\, \tilde Q_3 \left( \frac{w}{\e\, \cosh s}\right) ,
\]
where  $\tilde Q_2$ (resp. $\tilde Q_3$) fulfills properties similar to those enjoyed by $Q_2$ (resp. $Q_3$) in the statement of the result.

\medskip

The expansion of the normal vector $N^z_{eucl}(w)$ in terms of the function $w$ yields 
\[
N^z_{eucl} (w) = \tanh s + \frac{1}{\e\, \cosh^2 s}\,  \tilde L w + \tilde Q_2 \left( \frac{w}{\e \, \cosh s}\right) ,
\]
where $\tilde L$ (resp. $\tilde Q_2$) fulfills properties similar to those enjoyed by $L$ (resp. $Q_2$) in the statement of the result.

\medskip

Finally, we find  $z(w)$ in terms of the function $w$
\[
z(w)= 1 + \e \, s + \tanh s \, w  + \tilde Q_2' (w)  ,
\]
where $\tilde Q_{2}'$ is a nonlinear first order differential operator which fulfills properties similar to those enjoyed by $Q_2$ in the statement of the result.

\medskip

Writing $H_{hyp} = z \, H_{eucl} + N^z_{eucl}$, we obtain the desired expansion for $H_{hyp}$ in terms of $w$. \hfill $\Box$

\medskip

\begin{remark}
For notational convenience, we will write for short ${\cal L} w = Q (w)$ instead of (\ref{eq:nn}).
\end{remark}

\subsection{Mean curvature $1$ surfaces close to the truncated symmetrized catenoid}

We now fix $\mu \in (1,2)$. Recall that, for all $\e_{ij} \in (0,1)$ satisfying (\ref{eq:**}), we have  defined $s_{ij} >0$ by the identity
\[
\e_{ij} \, \cosh s_{ij}  = \rho .
\]
We set
\[
v_0 = \tanh s \int_0^s (\tanh t)^{-2}\int_0^t \tanh \zeta \, f_{ij} (\zeta)  \, d\zeta \, dt ,
\]
and we check directly that 
\[
\|v_0\|_{2, \alpha, 0}\leq c\, \e^2 \, \log^2 \e.
\]
Now, for all $\psi^i_\perp, \psi^j_\perp \in {\cal C}^{2, \alpha} (S^1)$ which are orthogonal to $1$, $e^{i\theta}$ and $e^{-i\theta}$  in the $L^2$ sense on $S^1$, we set ${\Psi^\perp_{ij}} : = (\psi^i_\perp, \psi_\perp^j)$ and define, for all $(s, \theta) \in [s_{ij}, s_{ij}]\times S^1$ 
\[
\tilde v_{\Psi^\perp_{ij}} (s) = {\cal P}\, \psi^i_\perp (s_{ij} -s)+ {\cal P}\, \psi^j_\perp (s+s_{ij}).
\]
Thanks to the result of Proposition~\ref{pr:7}, we get
\[
\|\tilde v_{\Psi^\perp_{ij}} \|_{2, \alpha, 2}\leq c \, \e^2 \, \|{\Psi^\perp_{ij}}\|_{2, \alpha}.
\]
Finally, we define $w_0 : = v_{\Psi^\perp_{ij}} + \tilde v_0$.

\medskip

Granted the above defined functions, we would like to solve the equation
\[
\left\{ 
\begin{array}{rllll}
{\cal L} (w_0+w)  & = &  Q(w_0+w) \qquad & \mbox{in}\qquad (-s_{ij}, s_{ij}) \times S^1\\[3mm]
w & \in &  \mbox{Span}\{ 1, e^{i\theta}, e^{-i\theta}\}\qquad &  \mbox{on}\qquad \{\pm s_{ij}\}\times S^1 .
\end{array}
\right.
\]
In order to do so, we use the result of Proposition~\ref{pr:6}, and we set
\[
\begin{array}{rllll}
{\cal S}_{\e_{ij}, \Psi^\perp_{ij}}  (w) & : = & \ds {\cal G}_{s_{ij}} \, \left( - {\cal L} v_{\Psi^\perp_{ij}} + \e_{ij}  \,  L(w_0+w)+ (\e_{ij} + \e_{ij}^2 \,  \cosh^2 s) \, Q_2\left(\frac{w_0 + w}{\e_{ij} \,  \cosh s}\right) \right. \\[3mm]
& & \qquad \ds \left. + \e_{ij} \, \cosh s \, Q_3 \left(\frac{w_0 + w}{\e_{ij} \, \cosh s }\right)\right) ,
\end{array}
\]
so that we can rewrite the above equation as a fixed point problem
\begin{equation}
w  = {\cal S}_{\e_{ij} , \Psi^\perp_{ij}} (w) ,
\label{eq:5-5}
\end{equation}
which we will solve by a standard contraction mapping argument. 

\medskip

To begin with, observe that, if we assume that $\|\Psi^\perp_{ij}\|_{{\cal C}^{2, \alpha}}\leq \kappa \, \e$, we have 
\[
\|{\cal L} \tilde v_{\Psi^\perp_{ij}} \|_{{\cal C}^{0, \alpha}_\mu}\leq c \, \e^3, \qquad \qquad \| \e_{ij} \, L (w_0)\|_{{\cal C}^{0, \alpha}_\mu}\leq c \, \e^3 \, \log^2 \e ,
\]
while 
\[
\left\| (\e_{ij}+ \e_{ij}^2 \, \cosh^2 s ) \, Q_2 \left( \frac{w_0}{\e_{ij} \cosh s}\right) \right\|_{{\cal C}^{0, \alpha}_\mu}\leq c \,  \e^3 \, \log^4 \e ,
\]
and 
\[
\left\| \e_{ij} \cosh s  \, Q_3 \left( \frac{w_0}{\e_{ij} \cosh s}\right) \right\|_{{\cal C}^{0, \alpha}_\mu}\leq c \,  \e^4 \, \log^6 \e ,
\]
where all constants do depend on $\kappa$. 

\medskip

It is then a simple exercise to show that, for all $\e$ small enough, there exists a unique solution of (\ref{eq:5-5}) which belongs to the ball of radius $c_0 \, \e^3 \, \log^4 \e$ in ${\cal C}^{2, \alpha}_\mu ([-s_{ij} , s_{ij}]\times S^1)$ and which is obtained as a fixed point for some contraction map. We leave the details to the reader.

\medskip

To summarize, we have produced a constant mean curvature $1$ surface which is close to the truncated symmetrized catenoid and which has two boundaries. Near those boundaries, this surface can be parameterized as a vertical graph over the $z=1$ plane. More precisely, close to the upper boundary, we can write this surface as the graph of 
\[
x \in B_{\rho}(a^j_\infty)-B_{\rho/2}(a^j_\infty)  \longrightarrow 1  - \e \, \log \e + \widehat W_{\psi^i_j} + \widehat V^i ,
\]
where $\widehat W_{\psi^i_j}$ is the harmonic extension of the function 
\[
\psi^i_j : = d^i_j - (\e_{ij} -\e) \, \log \e -\frac{2}{|\xi^i_j|^2} (x- a^j_\infty) \cdot \xi^i_j  + \psi^i_\perp ( \frac{x-a^j_\infty}{\rho})  ,
\]
in $B_{\rho} (a^j_\infty)$. Recall that we have set $\xi^i_j :=  a^i_j-a^j_i$. 

\medskip

Observe that, after an inversion centered at $a^j_i$, we obtain the parameterization of the lower end as
\[
x \in B_{\rho}(a^j_\infty)-B_{\rho/2}(a^j_\infty)  \longrightarrow 1 - \e \, \log \e + \widehat W_{\psi^j_i} + \widehat V^j,
\]
where $\widehat W_{\psi^j_i}$ is the harmonic extension of the function 
\[
\psi^j_i : = - d^i_j - (\e_{ij} -\e) \, \log \e - \frac{2}{| \hat \xi^i_j|^2} (x- x^j_i(\e)) \cdot  \hat \xi^i_j  + \psi^j_\perp (\frac{x-x^i_j(\e)}{\rho})  ,
\]
where this time
\[
\hat \xi_{ji}:= \hat a^j_i - \hat a_j^i ,
\]
with $\hat a^i_j$ and $\hat a^j_i$ given by
\[
\hat a_j^i := \frac{a^i_j -a_\infty^j}{|a^i_j -a_\infty^j|^2} \qquad \mbox{and} \qquad \hat a_i^j := \frac{a^j_i -a_\infty^j}{|a^j_i -a_\infty^j|^2}.
\]
The key point is that this formula involves $-d^i_j$ and not $+d^i_j$ as in the previous one, hence, by adjusting properly the parameters $\e_{ij}$, $d^i_j$, $a^i_j$ and $a_i^j$ we can prescribe any boundary data  the lower or upper boundary provided it is bounded by $\kappa \, \e$ in ${\cal C}^{2, \alpha}(S^1)$ norm. 

\medskip

In the above formula the mappings $V^i$ and $V^j$ depend smoothly on the data $\psi^i_j, \psi_i^j$ and  are bounded by a constant (independent of $\kappa$) times $\e$ in ${\cal C}^{2, \alpha} (B_{\rho}-B_{\rho / 2})$ topology. Moreover, we have
\begin{equation}
\ds \|  \widehat V^i_{\psi^i_j, \psi_i^j} -  \widehat V_{\tilde \psi^i_j, \tilde \psi_i^j}\|_{{\cal C}^{2, \alpha}(B_\rho- B_{\rho/2})} \leq \ds c \, \frac{1}{|\log \e|}\, ( \| \psi^i_j - \tilde \psi^i_j\|_{{\cal C}^{2, \alpha}} + \| \psi^j_i - \tilde \psi^j_i\|_{{\cal C}^{2, \alpha}}  ) .
\label{eq:cont-2}
\end{equation}

\medskip

We will denote by $C_{ij} (\psi^i_j, \psi^j_i)$ the surface constructed in this section. 

\section{The gluing procedure}

We will fix $\kappa >0$ large enough at the end of this section and apply the results of the previous sections. In  particular, $\delta$ is fixed in $(-1, 0)$ and $\mu$ is fixed in $(1,2)$. There exists $\e_0 >0$ and for any collection of boundary data $\phi^i_j, \psi^i_j \in {\cal C}^{2, \alpha}(S^1)$ satisfying $\|\phi^i_j\|_{2, \alpha}\leq \kappa \, \e$ and $\|\psi^i_j\|_{2, \alpha}\leq \kappa \, \e$, we can find {\small CMC}-$1$ surfaces $M_i ( \Phi_i)$ which are close to the horosphere ${\cal H}^i$ and {\small CMC}-$1$ surfaces $C_{ij} (\psi^i_j , \psi^j_i)$ which are close to truncated symmetrize catenoids. Our aim will now be to find $\phi^i_j$ and $\psi^i_j$  in such a way that 
\[
\left( \cup_i M_i ( \Phi_i) \right) \cup \left( \cup_{i,j} C_{ij} (\psi^i_j , \psi^j_i)\right) ,
\]
is a ${\cal C}^{1}$ surface whose mean curvature is constant away from the boundaries of the different pieces. Then, we can  apply standard regularity theory to show that this surface is in fact ${\cal C}^\infty$ since it is ${\cal C}^1$ and has mean curvature equal to $1$. 
 
\medskip

By construction, the two surfaces  $M_i ( \Phi_i) $ and $C_{ij} (\psi^i , \psi^j)$ are graphs over the $z=1$ plane near their boundary which is close to the point $p_{ij}$. More precisely, $M_i ( \Phi_i)$ is the graph of 
\[
x \in B_{2 \rho} (a^j_\infty) - B_{\rho}(a^j_\infty) \longrightarrow  1 - \e\, \log \e +  W_{\Phi_i} + V_{\Phi_i},
\]
and the upper boundary of $C_{ij} (\psi_j^i, \psi_i^j)$ is the graph of 
\[
x \in B_{\rho}(a^j_\infty) - B_{\rho/2}(a^j_\infty) \longrightarrow  1  - \e \, \log \e + \widehat W_{\psi^i_j} + \widehat V^i_{\psi^i_j, \psi^j_i}.
\]

Hence, to produce a ${\cal C}^1$ surface, it remains to ask that the Dirichlet data of these two graphs coincide, so that it already ensures that the surface is ${\cal C}^0$ and also that the Neumann data of these two graphs coincide, so that it ensures that the surface will be of class ${\cal C}^1$.

\medskip

Let us denote by $n_i$ the cardinal of $I_i$. We claim that the mapping ${\cal U}$
\[
{\cal U} :  (\phi^j)_j \in ({\cal C}^{2, \alpha}(S^1 ))^{n_i} \longrightarrow \rho \, (\del_r (W_{\Phi_j} - \widehat W_{\phi^j})(\rho \, \cdot +a^{j'}_\infty ) )_{j'} \in ({\cal C}^{1, \alpha} (S^1))^{n_i} ,
\]
is an isomorphism. Indeed, this mapping is a linear first order elliptic pseudo-differential operator with principal symbol $ - 2 \, |\xi| $. Therefore, in order to check that it is an isomorphism, it is enough to prove that it is injective. Now if we assume that ${\cal U} ((\phi^j)_j)=0$ then the function $w$ defined by $w := \widehat W_{\Phi_j}$ in $\Omega_{i}$ and $w: =  \widehat W_{\phi^j}$ in $B_\rho (a^j_\infty)$ is a solution of $\Delta w =0$ in ${\R}^2$, and furthermore, $w$ belongs to ${\cal C}^{2, \alpha}_{\delta} ({\R}^2) \oplus {\cal K}_0$. As already mentioned in Remark~\ref{re:inj}, this implies that $w \equiv 0$ and, as a consequence, $h\equiv 0$.  

\medskip

Using the above claim,  it is easy to see that the problem reduces to a fixed point problem 
\[
(\phi_j^i, \psi^i_j )_{ij} = {\bf C}_\e ( (\phi_j^i, \psi^i_j )_{ij}) ,
\]
in $F : = ({\cal C}^{2, \alpha}(S^1))^N$, where $N = 2 \, \sum_i n_i$. However,  (\ref{eq:cont-1}) and (\ref{eq:cont-2}) imply that, the constant $\kappa$ being fixed large enough, the mapping ${\bf C}_\e :  F \longrightarrow F$ is a contraction mapping defined in the ball of radius $\kappa \, \e$ of $F$ into itself, for all $\e$ small enough. Hence, we have obtained a fixed point of the mapping ${\bf C}_\e$. This completes the proof of both Theorem~1 and  Theorem~3.

\section{The moduli space theory and the nondegeneracy of the solutions}

Assume that $M$ is an orientable {\small CMC}-$1$ surface all of whose ends  $E_1, \ldots , E_n$ are asymptotic to the end of a catenoid cousin. Up to some isometry,  each end $E_i$ can be parameterized as a graph of a function of the form
\[
{\R}^2 - B_\rho \ni x \longrightarrow  r^{2 - t_i} + {\cal O}(r^{4-3t_i}) ,
\]
where $t_i \in (1, +\infty)$. As we have shown in section \S 3, there are, for each end, $6$ linearly independent Jacobi fields which correspond to the $6$ different geometric transformations. We shall denote them by 
\[
h^{j, \pm}_i \qquad \mbox{for}\qquad j=-1, 0, +1 \qquad \mbox{and}\qquad i=1, \ldots, n ,
\]
where $h^{j, \pm}_i \sim h^{j, \pm}$, at $\infty$.
 
\medskip

We now decompose $M$ into slightly overlapping pieces which are a compact piece $M^c$ and the ends $E_i$. Furthermore, we ask that, for each $i=1, \ldots , n$, the set $M^c \cap E_i$ is diffeomorphic to $[0,1]\times S^1$. With this decomposition, we give the~:
\begin{definition}
The function space ${\cal E}^{k,\alpha}_\delta  (M)$ is defined to be the space of all functions $w\in {\cal C}^{k, \alpha}_{loc}(M)$ for which  the following norm is finite
\[
|w |_{{\cal E}^{k, \alpha}_\delta} : = \sum_{i=1}^n \|w_{|_{E_i}}\|_{{\cal C}^{k, \alpha}_\delta ({\R}^2 - B_\rho) }+ [w_{|_{M^c}}]_{k, \alpha, M^c} ,
\]
where $\|\,\, \|_{{\cal C}^{k, \alpha}_\delta }$ is the norm defined in Definition~\ref{de:4}. 
\end{definition} 

 \medskip

{\bf Moduli space theory :} To begin with, for any {\small CMC}-$1$ surface $M$, we denote by ${\cal L}_{M}$ the Jacobi operator about $M$, that is the linearized mean curvature operator about $M$ with respect to the (hyperbolic) normal vector field. Observe that, if all the ends of $M$ are asymptotic to ends of catenoid cousins, the set of indicial roots of ${\cal L}_{M}$  is given by ${\mathbb Z}$. Recall that these indicial roots determine the asymptotic behavior of any solution of the homogeneous problem ${\cal L}_M w=0$ near the ends. Following the analysis of \S 4, we can show that, for all  $\delta \notin {\mathbb Z}$,  the operator ${\cal L}_M : {\cal E}^{2, \alpha}_\delta (M) \longrightarrow {\cal E}^{0, \alpha}_{\delta -2} (M)$ is injective if and only if the operator ${\cal L}_M :{\cal E}^{2, \alpha}_{-\delta} (M) \longrightarrow {\cal E}^{0, \alpha}_{-\delta -2} (M)$ is surjective.

\medskip 

We can now give the precise definition of nondegeneracy.
\begin{definition}
We will say that a {\small CMC}-$1$ surface $M$ is unmarked nondegenerate, if the linearized mean curvature operator with respect to the (hyperbolic) normal vector field
\[
{\cal L}_M : {\cal E}^{2, \alpha}_\delta (M)\longrightarrow {\cal E}^{0, \alpha}_\delta(M) ,
\]
is injective for all $\delta \in (-\infty, -1)$.  We will say that a {\small CMC}-$1$ surface $M$ is marked nondegenerate, if 
\[
{\cal L}_M : {\cal E}^{2, \alpha}_\delta (M)\longrightarrow {\cal E}^{0, \alpha}_\delta(M) ,
\]
is injective for all $\delta \in (-\infty, 0)$. 
\end{definition}
For example, any horosphere is both marked and unmarked nondegenerate.

\medskip

Let us assume that $M$ is a {\small CMC}-$1$ surface which is marked nondegenerate, in which case we set $\iota =0$, or unmarked nondegenerate, in which case we set $\iota =1$. The deficiency spaces are defined by 
\[
{\cal D}_0 := \ds \oplus_{i=1, \ldots,n} \mbox{Span} \{ \, \chi \, h^{0, \pm}_i \},
 \]
and 
\[
{\cal D}_1 : =  \ds \oplus_{i=1, \ldots,n} \mbox{Span} \{ \, \chi \, h^{j, \pm}_i \, : \, j=-1, 0, 1\} ,
\]
where, as usual, $\chi$ is a cutoff function equal to $0$ in some sufficiently large ball and equal to $1$ outside a larger ball. We fix $\delta \in (-1- \iota, -\iota)$ and we define  ${\cal N}_\iota$  to be the trace of the kernel of ${\cal L}_M :{\cal E}^{2, \alpha}_{- \delta} (M) \longrightarrow {\cal E}^{0, \alpha}_{-\delta -2} (M)$ in ${\cal D}_\iota$, that is ${\cal N}_\iota$ is a $n \, (1+2 \iota )$ dimensional subspace of ${\cal D}_\iota$ such that
\[
\mbox{Ker} \, {\cal L}_{M} \subset   {\cal E}^{2, \alpha}_{\delta} (M) \oplus {\cal N}_\iota .
\]
Finally, we define ${\cal K}_\iota$ to be a $n\, (1+2\iota)$ dimensional subspace of ${\cal D}_\iota$ such that 
\[
{\cal D}_\iota ={\cal N}_\iota \oplus {\cal K}_\iota,
\]
then, we can show that 
\[  
{\cal L}_M : {\cal E}^{2, \alpha}_\delta (M) \oplus {\cal K}_\iota \longrightarrow {\cal E}^{0, \alpha}_\delta (M),
\]
is an isomorphism. Using this, the results of Theorem~\ref{th:2} and Theorem~\ref{th:4} follow at once by modifying the analysis of \cite{Kus-Maz-Pol}, or the analysis of \S 4.1.  The proof being almost identical to what is done in \cite{Kus-Maz-Pol}, we leave the details to the reader and concentrate on the proof of the nondegeneracy of the solutions constructed. 

 \medskip

{\bf Unmarked nondegeneracy :} We prove that the {\small CMC}-$1$ surfaces of genus $0$ we have obtained in this paper are marked nondegenerate, for all $\e$ small enough. In particular, this will imply that these surfaces are regular points of their respective moduli space. The proof is by contradiction.  Assume that, for a sequence $\e_k$ tending to $0$, the operator ${\cal L}_{M_{\e_k}}$ is not injective on ${\cal E}^{2, \alpha}_\delta (M_{\e_k})$, for some $\delta \in (-2,-1)$. If this is so, there exists for each $k$ some nontrivial function $w_k\in  {\cal E}^{2, \alpha}_\delta (M_{\e_k})$ such that ${\cal L}_{M_{\e_k}} w_k= 0$. 

 \medskip

By construction, we may decompose $M_{\e_k}$ into the union of a  pieces $H^k_i$ which are small normal geodesic graphs over the horospheres and pieces $C^k_{ij}$ which are small normal geodesic graphs over truncated, rescaled symmetrized catenoids. We define on each $M_{\e_k}$ some weight function $\gamma_k >0 $, as follows~:
\begin{itemize}
\item $\gamma_k \sim (1 + r^2)^{\delta /2}$ on $H^k_i$,
\item $\gamma_k \sim (\e_k  \, \cosh s)^\delta $ in $C^k_{ij}$,
\end{itemize}
where $f \sim g$ means that $1/2 \leq f/g\leq 2$. We have taken obvious coordinates to parameterize the surfaces $H^k_i$  which is close to  horospheres (which, up to an isometry,  can be assumed to be the horizontal plane $z=1$) and the surfaces $C^k_{ij}$ which are close to a truncated symmetrized catenoid (which up to an isometry can be assumed to be vertical). 

\medskip

Having defined this weight function, we normalize the sequence $w_k$ so that 
\[
\sup_{M_{\e_k}} \, \gamma_k^{-1}\, w_k =1.
\]
As already mentioned, the set of indicial roots of ${\cal L}_{M_{\e_k}}$ at each end is given by ${\mathbb Z}$. Hence any bounded solution of ${\cal L}_{M_{\e_k}}w=0$ which belongs to the space ${\cal E}^{2, \alpha}_\delta (M_{\e_k})$ decays like a constant times $r^{-2}$ at each end. This implies that the above supremum is achieved (say at some point $p_k \in M_{\e_k}$). We now distinguish a few cases according to the behavior of the sequence $p_k$. By construction, as $k$ tends to $\infty$, the sequence of surfaces $M_{\e_k}$ converges to the union of horospheres ${\cal H}^i$ which are connected at points $p_{ij}$. Since we are working in the case where the genus of $M_{\e_k}$ is $0$, we can label the horospheres in such a way that ${\cal H}^i$ is linked to ${\cal H}^{i+1}$ for all $i=1, \ldots,  n-1$.

\medskip

\noindent
{\em Case 1 :} Assume that, up to a subsequence, the sequence $p_k$ converges to some point $p_\infty$ which is the end of one of the horospheres, say ${\cal H}^i$. For $k$ large enough $p_k$ corresponds to some point of coordinate $x_k \in{\R}^2$, if we assume that the limit horophere ${\cal H}^i$ is the horizontal plane $z=1$. Extracting some subsequences, if this is necessary, we find that the sequence of rescaled functions
\[
\tilde w_k  := |x_k|^{-\delta} \, w_k (|x_k|\, \cdot) ,
\]
converges uniformly on any compact of ${\R}^2\setminus \{0\}$ to a nontrivial solution of 
\[
\Delta w_\infty =0,
\]
in ${\R}^2 -\{0\}$. Moreover, $w_\infty$ is bounded by a constant times $r^\delta$. But this is easily seen to be impossible by considering the Fourier decomposition of $w_\infty$ in the $\theta$ variable, since we have assume $\delta \notin {\mathbb Z}$.

\medskip

\noindent
{\em Case 2 :} Assume that, up to a subsequence, the sequence $p_k$ converges to some point $p_\infty \in {\cal H}^i$ which is not a point where the horospheres are connected. For $k$ large enough $p_k$ corresponds to some point of coordinate $x_k \in{\R}^2$, if we assume that the limit horophere ${\cal H}^i$ is the horizontal plane $z=1$. Extracting some subsequences, if this is necessary, we may assume that the sequence $w_k$ converges uniformly on any compact of the horosphere ${\cal H}^i \setminus \{p^i_{i+1}, p^i_{i-1}\}$, if $i\in \{2, \ldots, n-1 \}$  (or of ${\cal H}^1 \setminus \{p^1_{2}\}$ if $i=1$, or of ${\cal H}^n \setminus \{p^n_{n-1}\}$ if $i=n$)   to a solution of 
\[
\Delta w^i_\infty =0 ,
\]
in  ${\cal H}^i \setminus \{p^i_{i+1}, p^i_{i-1}\}$. Furthermore, in the collection $(w^1_\infty, \ldots, w_\infty^n)$, at least one of the functions is nontrivial. Observe that we have identified each horosphere with the hyperplane $z=1$, which we might well after having performed an isometry. Moreover, $w^i_\infty$ is bounded by a constant times $r^\delta$ at $\infty$ and is bounded by a constant times $|x-x^i_{i+1}|^{-\delta}$ near $x^i_{i+1}$ and by a constant times  $|x-x^i_{i -1}|^{-\delta}$ near $x^i_{i-1}$, where $x^i_{i+1}$ is the projection of $p^i_{i+1}$ over $\del {\HH}^3$.

\medskip

As in the previous case, it is easy to see that necessarily $w^1_\infty =0$ and also that $w^n_\infty =0$. We shall now prove, as in \cite{Maz-Pac-Pol-1}, that all $w^i_\infty \equiv 0$. For example, $w^2_\infty$ satisfies 
\[
\Delta w^2_\infty = a_1 \, \delta_{x^2_1} + a_3 \, \delta_{x^2_3} + b_1\cdot \nabla \delta_{x^2_1} + b_3\cdot \nabla \delta_{x^2_3} ,
\]
for some $a_1, a_3 \in {\R}$ and some $b_1, b_3 \in {\R}^2$. To begin with let us show that $a_1=0$ and also that $b_1=0$, it will follow immediately, as in case 1, that $w^2_\infty \equiv 0$. Then, a simple induction will show that all $w^i_\infty \equiv 0$, which is the desired contradiction. 

\medskip

To show that $a_1=0$ and $b_1 =0$, we use the fact that, in any annular neighborhood of $x^2_1$, the function $w^2_\infty$ is the limit of the functions $w_k$ which are defined on $M_{\e_k}$. Now, up to an isometry, we can assume that ${\cal H}^2$ is the horizontal plane $z=1$, ${\cal H}^1$ is the horosphere of radius $1/2$ centered at $(0,1/2)$ and that  $x^2_1=  0$. In particular over an annular region around $x^2_1$ which is close to $p_{12}$, the surfaces $M_{\e_k}$ are multiple valued graphs over an annulus $B_{2\rho} - B_\rho$, one of the graphs corresponding to the piece of $M_{\e_k}$ which is close to ${\cal H}^1$ and the other corresponding to the piece of $M_{\e_k}$ which is close to ${\cal H}^2$.  Letting the group of isometries act of the surface $M_{\e_k}$, we can explicitely describe some Jacobi fields.  Of interest will be the Jacobi field $\tilde h_k^{-,0}$  which corresponds to a dilation centered at the origin and the Jacobi fields $\tilde h^{1, \pm}_k$ which correspond to the isometry which keep the point $(0,1) \in \del {\HH}^3$ fixed and send $\infty $ to some point $a \in \del {\HH}^3$.  The only information we will need is that 
\[
\tilde h^{-,0}_k \sim 1,
\]
\[
\tilde h^{1,+}_k \sim r \, \cos \theta  \qquad \mbox{and}\qquad  \tilde h^{1,-}_k \sim r \, \sin \theta ,
\]
in the annular region of $M_{\e_k}$ which is close to ${\cal H}^2$.

\medskip

Now consider 
\[
h \, ({\cal L}_{M_{\e_k}} w_k ) -  w_k \, ( {\cal L}_{M_{\e_k}}  h ) =0, 
\]
where $h$ is any one of the above described Jacobi fields and integrate this by parts over the catenoidal neck between the two annular regions. We take advantage of the fact that we already know that $w_k$ converges to $0$ on the annular region close to ${\cal H}^1$ and also that $w_k$ converges to $w^2_\infty$ in the anular region close to ${\cal H}^2$. Moreover, we know that, near $x^2_1$ we have the expansion
\[
w^2_\infty =  - \frac{a_1}{2\, \pi}\, \log r  - \frac{b_1}{2\, \pi}\cdot \nabla \log r + {\cal O}(1).
\]
Taking $h=\tilde h^{-,0}_k$ and letting $k$ tend to $\infty$, we find that $a_1=0$, finally, taking $h= \tilde h_k^{1, \pm}$ and letting $k$ tend to $\infty$ we find that $b_1=0$. The reader can find the details of a similar argument in \cite{Maz-Pac-Pol-1}.

\medskip

\noindent
{\em Case 3 :} Assume that, up to a subsequence, the sequence $p_k$ tends to a point where two horospheres are connected. Hence, for $k$ large enough, the point $p_k$ corresponds to a point $(s_k, \theta_k)$ in the parameterization of the connecting catenoid.  In this last case still distinguish two subcases according to whether $s_k$ remains bounded or tends to $\pm \infty$. If $s_k$ tends to $-\infty$, we define the sequence of rescaled functions 
\[
\tilde{w}_i (s, \theta) : = e^{\delta \, s_k}\, w_i( s + s_k, \theta),
\]
which, up to a subsequence, may assumed to converge to a nontrivial solution of 
\[
(\del^2_s + \del_\theta^2) \, w_\infty =0\qquad \mbox{on}\qquad {\R}\times S^{n-1},
\] 
which  is bounded by a constant times $e^{- \delta s}$. Expending $w_\infty$ in Fourier series, one easily checks that this is not possible. A similar argument holds when $s_k$ tends to $+\infty$. Finally, if $s_k$ converges to $s_* \in {\R}$, we simply get that $\e_k^{-\delta} \, w_k$ converges to a nontrivial solution of 
\[
{\cal L} \, w_\infty =0 \qquad \mbox{on}\qquad {\R}\times S^{n-1},
\]
and is bounded by a constant times $(\cosh s)^{\delta}$. Again,  to see that this is not possible we decompose $w_\infty$ as $\sum_m w_m \, e^{im \theta}$. Thanks to the choice of $\delta$,  we see that $w_0=w_{\pm 1} =0$ since all nontrivial solutions of the homogeneous problem ${\cal L} w=0$ corresponding to the eigenspace $m=0, \pm 1$ decay at most like $\cosh^{-1}s$ at $\infty$, hence can't be bounded by $(\cosh s)^{\delta}$ for some $\delta \in (-2, -1)$. Now, when restricted to the eigenspaces corresponding to $m \geq 2$,  as already mentioned in the proof of Proposition~\ref{pr:4}, the operator ${\cal L}$ satisfies the maximum principle and it is easy to show that $w_m=0$ for all $m \in {\mathbb Z}$.  Again a contradiction.

\medskip

\noindent
Since we have ruled out every possible case, the proof of fact that the surfaces constructed in this paper are unmarked nondegenerate, provided $\e$ is small enough, is complete.

\medskip

{\bf Generic marked nondegeneracy :}  To finish, we prove that the solutions we construct are  marked nondegenerate for small values of $\e$ and for any generic configuration of the ends. The proof of this fact is borrowed from \cite{Maz-Pac-0}. Indeed, as discussed in \cite{Maz-Pol-Uhl}, near smooth points $M$ of the unmarked moduli space ${\cal M}^u_{g,n}$, there is a real analytic fibration $\pi: {\cal M}^u_{g,n} \rightarrow {\cal C}_n$ onto the configuration space of $n$ distinct points in $\del {\HH}^3$ corresponding to the ends of the surface. Sard's theorem implies the surjectivity of the differential of $\pi$ for all points in a generic fiber. Surjectivity of this differential (which we assume from now on) is equivalent to the fact that, in the unmarked moduli space one can consider  the ends of the surface to be independent parameters, at least locally at any point $M$ in the fiber. In other words one can freely move any of the ends in $\del {\HH}^3$. Differentiation with respect to any of these $2n$ parameters gives rise to a Jacobi field which is bounded by a constant time $- \log r$ near each end except one where is behaves like  $h^{+1, \pm }$. We will denote by $g^{+ 1, \pm}_k$ the Jacobi fields obtained this way, the indices referring to the fact that this is the Jacobi field which behaves like  $h^{+1, \pm}$ at the $k$-th end. 

\medskip

Now, if $w$ is any solution of the homogeneous equation ${\cal L}_M w=0$ in $M$, which belongs to ${\cal E}^{2, \alpha}_\delta$ for some $\delta <0$.  Inspection of the asymptotic behavior of $w$ near each end shows that $w \in {\cal E}^{2, \alpha}_{-1}$. Furthermore, near the end $E_k$ we have the following expansion
\[
w  = a_k^+ \,  h^{-1, +}_k + a^-_k \, h^{-1,-}_k + {\cal O} (r^{-2}) ,
\]
where $h^{-1,\pm}_k$ have been defined at the beginning of this section.  For fixed $k=1, \ldots, n$, we can integrate $w \, ({\cal L}_M g^{+1, \pm}_k) - g^{+1, \pm}_k \, ({\cal L}_M w ) = 0$ over a sequence of compacts which exhaust $M$. Passing to the limit, we obtain $a_k^\pm=0$. Hence $w \in {\cal E}^{2, \alpha}_{\delta'}$ for some $\delta' <-1$. Since the surface $M$ has been assumed to be unmarked nondegenerate, we conclude that $w \equiv 0$. This finishes the proof that the surfaces are marked nondegeneracy of all points in the fiber.


\begin{thebibliography}{99}
\bibitem{Tou} J.L. Barbosa and R. Sa Earp. {\em Prescribed mean curvature hypersurfaces in ${\mathbb H}^{n+1}$ with convex planar boundary II}. S{\'e}minaire de th{\'e}orie spectrale et g{\'e}om{\'e}trie de Grenoble. {\bf 16}, (1998), 43-79.

\bibitem{Ber-Ros} J. Berglund and W. Rossman. {\em Minimal surfaces with catenoidal ends}. Pacific J. Math. {\bf 171}, 2, (1985),  353-371.

\bibitem{Bry} R.L. Bryant. {\em Surfaces of mean curvature one in hyperbolic space}. Th{\'e}orie des vari{\'e}t{\'e}s minimales et applications, Ast{\'e}risque 154-155, Soc. Math. France, Paris (1987).

\bibitem{Col-Hau-Ros} P. Collin, L. Hauswirth and H. Rosenberg. {\em The geometry of finite topology surfaces properly embedded in hyperbolic space with constant mean curvature one}, preprint.

\bibitem{Col-Ros} P. Collin and H. Rosenberg. In preparation.

\bibitem{Cos-Ros} C.P. Cosin and A. Ros. {\em A Plateau problem at infinity for properly immersed minimal surfaces with finite total curvature}, Preprint, 1998.

\bibitem{Del} C. Delaunay. {\em Sur la surface de revolution dont la courbure moyenne est constante}.  J. de Math{\'e}matiques, {\bf 6} (1841) 309-320.

\bibitem{Die-Hil} U. Dierkes, S. Hildebrandt, A. K{\"u}ster and O. Wohlrab. {\em Minimal surfaces I. Boundary value problems}. Springer-Verlag 295 (1992).

\bibitem{Fak-Pac} S. Fakhi and F. Pacard. {\em Existence of complete minimal hypersurfaces with finite total curvature}. Manuscripta Mathematica. {\bf 103},  (2000), 465-512.

\bibitem{KGB} K. Grosse-Brauchmann. {\em New surfaces of constant mean curvature}. Math. Zeit. {\bf 214}, (1993), 527-565.

\bibitem{Gro-Kus-Sul} K. Grosse-Brauchmann, R.B. Kusner and J.M. Sullivan. {\em Classification of embedded constant mean curvature surfaces with genus zero and three ends},  preprint (2000).

\bibitem{Jor-Mee} L. Jorge and W.H. Meeks III. {\em The topology of complete minimal surfaces of finite total curvature}. Topology, {\bf 2}, (1983), 203-221.

\bibitem{Kap} N. Kapouleas. {\em Complete constant mean curvature surfaces in Euclidean three space}. Ann. of Math. (2), {\bf 131}, 2, (1990), 239-330.

\bibitem{Kat-Ume-Yam} S. Kato, M. Umehara and K. Yamada. {\em An inverse proble of the flux for minimal surfaces}. Indiana Univ. Math. J. {\bf 46}, (1997), 529-559.

\bibitem{Kor-Kus-Mee-Sol} N.J. Korevaar, R.B. Kusner, W.H. Meeks and B. Solomon. {\em Constant mean curvature surfaces in hyperbolic space}.  Amer. J. Math. {\bf 114}, 1, (1992), 1-43.

\bibitem{Kus-Maz-Pol} R.B. Kusner, R. Mazzeo and D. Pollack. {\em The moduli space of complete embedded constant mean curvature surfaces}. Geom. Funct. Anal. {\bf 6}, (1996), 120-137.

\bibitem{Lev-Ros} G. Levitt and H. Rosenberg. {\em Symmetry of constant mean curvature hypersurfaces in hyperbolic space}. Duke Math. J. {\bf 52}, 1, (1985), 53-59.

\bibitem{Maz-Pac-0} R. Mazzeo and F. Pacard. {\em Constant scalar curvature metrics with isolated singularities.} Duke Math. J, {\bf 99}, (1999), 3, 353-418.

\bibitem{Maz-Pac} R. Mazzeo and F. Pacard. {\em Constant mean curvature surfaces with Delaunay ends}.  Comm. Analysis and Geometry. {\bf 9}, 1, (2001), 169-237.
 
\bibitem{Maz-Pac-Pol-1} R. Mazzeo, F. Pacard and D. Pollack. {\em Connected sums of constant mean curvature surfaces in Euclidean 3 space}. To appear in J. Reine Angew. Math. (2001).

\bibitem{Maz-Pac-Pol-Rat} R. Mazzeo, F. Pacard, D. Pollack and J. Ratzkin. In preparation.

\bibitem{Maz-Pol} R. Mazzeo and D. Pollack. {\em Gluing and moduli for some noncompact geometric problems}. Geometric Theory of Singular Phenomena in
Partial Differential Equations, Symposia Mathematica {\bf 38},   Cambridge Univ. Press (1998) 17-51. 

\bibitem{Maz-Pol-Uhl} R. Mazzeo, D. Pollack and K. Uhlenbeck. {\em Moduli spaces of singular Yamabe metrics}. J. Amer. Math. Soc. {\bf 9} (1996) 303-344.

\bibitem{Oss} R. Osserman.  {\em A survey of minimal surfaces}. Dover Publications, New York, 2nd Edition, (1986).

\bibitem{Pac-Riv} F. Pacard and T. Riviere. {\em Linear and nonlinear aspects of vortices : the Ginzburg-Landau model}. Progress in Nonlinear Differential Equations, 39, Birk{\"a}user. 342 pp. (2000).

\bibitem{Per-Ros} J. Perez and A. Ros. {\em The space of properly embedded minimal surfaces with finite total curvature}. Indiana Univ. Math. J. {\bf 45},  (1996), 177-204.

\bibitem{Ros-Sat} W. Rossman and K. Sato.  {\em Constant mean curvature surfaces with two ends in hyperbolic space}. Experimental Math.  {\bf 7}, 2, (1998), 101-119.

\bibitem{Ros-Ume-Yam} W. Rossman, M. Umehada and K. Yamada. {\em Irreducible constant mean curvature $1$ surfaces in hyperbolic space with positive genus}. Tohoku Math. J. {\bf 49}, (1997), 449-484. 

\bibitem{Ros-Ume-Yam-2} W. Rossman, M. Umehada and K. Yamada. {\em Mean curvature $1$ surfaces in hperbolic $3$-space with low total curvature I}. arXiv:math.DE/0008015

\bibitem{Spi} M. Spivak. {\em A comprehensive introduction to differential geometry}. Second edition, Publish or Perish, Inc. (1979).

\bibitem{Ume-Yam} M. Umehara and K. Yamada. {\em Complete surfaces of constant mean curvature-$1$ in the hyperbolic 3-space}. Ann. of  Maths. {\bf 137}, (1993), 611-638.

\bibitem{Yan} S.D. Yang. {\em A connected sum construction for complete minimal surfaces with fnite total curvature}. Comm. Analysis and Geometry. {\bf 9}, 1, (2001), 116-168. 
\end{thebibliography}
\end{document}